\documentstyle[11pt, epsfig, pb-diagram]{article}
\input amssym.def
\title{\bf Parallel spinors and connections with skew-symmetric torsion in string theory
\footnote{Supported by the SFB 288 of the DFG. S.I. is a member of the EDGE,Research Training Network HPRN-CT-2000-00101, supported by the European Human Potential Programme.S.I. thanks the Humboldt University Berlin and ICTP, Trieste
for the support and excellent enviroments.}}
\author{Thomas Friedrich (Berlin) and Stefan Ivanov (Sofia)}
\date{\today}

\begin{document}

%==========================
\newcommand{\D}{\displaystyle}
\newcommand{\upsp}{\phantom{l}}
\newcommand{\downsp}{\phantom{q}}
%==========================
\newfont{\graf}{eufm10}
\newcommand{\son}{\mbox{\graf so}(n)}
\newcommand{\altg}{\mbox{\graf g}}
\newcommand{\altm}{\mbox{\graf m}}
\newcommand{\bdm}{\begin{displaymath}}
\newcommand{\edm}{\end{displaymath}}

\def\haken{\mathbin{\hbox to 6pt{%
                 \vrule height0.4pt width5pt depth0pt
                 \kern-.4pt
                 \vrule height6pt width0.4pt depth0pt\hss}}}
    \let \hook\intprod

\maketitle

\mbox{} \hrulefill \mbox{}\\

\newcommand{\vol}{\mbox{vol} \, }
\newcommand{\grad}{\mbox{grad} \, }

\begin{abstract} We describe all almost contact metric, almost hermitian
and $G_2$-structures admitting a connection with totally skew-symmetric
torsion tensor, and prove that there exists at most one such connection.
We investigate its torsion form, its Ricci tensor, the Dirac operator
and the $\nabla$-parallel spinors. In particular, we obtain solutions of
the type $II$ string equations in dimension $n=5,6$ and $7$.  
\end{abstract}

{\small
{\it Subj. Class.:} Special Riemannian manifolds, Spin geometry, String theory
} \\

\setcounter{section}{0}

\mbox{} \hrulefill \mbox{}\\
%----------------------------------------------------------------------------
\section{Introduction}
%----------------------------------------------------------------------------
Linear connections preserving a Riemannian metric with totally 
skew-symmetric torsion recently became a subject of interest in theoretical and mathematical physics. For example, the target space of supersymmetric sigma 
models with Wess-Zumino term carries a geometry of a metric connection with 
skew-symmetric torsion \cite{GHR,HP1,HP2} (see also \cite{OP1} and references 
therein). In supergravity theories, the  geometry of the moduli space of a 
class of black holes is carried out by a metric connection with  skew-symmetric torsion \cite{GPS}. The geometry of NS-5 brane solutions of type II supergravity theories is generated by a metric connection with skew-symmetric torsion 
\cite{PT1,PT2,P1}. The existence of parallel spinors with respect to a metric 
connection with skew-symmetric torsion on a Riemannian spin manifold is of  
importance in string theory, since they are associated with some string 
solitons (BPS solitons)  \cite{P1}. Supergravity solutions that preserve some 
of the supersymmetry of the underlying theory have found many applications in 
the exploration of perturbative and non-perturbative properties of string 
theory. An important example is the AdS/CFT correspondence, also known  as 
the Maldacena conjecture, which conjectures duality between the supersymmetric 
background and a certain superconformal field theory \cite{Mal,GKP,W}. \\

In type II string  theory one investigates manifolds
$N^k\times M^{10-k}$, where $N^k$ is a $k$-dimensional space-time and 
$M^{10-k}$ is a Riemannian manifold equipped with an additional 
structure. Indeed, the basic model is a $5$-tuple $(M^n,g,H,\Phi,\Psi)$, 
where $g$ is a Riemannian metric, $H$ is a $3$-form, $\Phi$ is the socalled 
dilation function, and $\Psi$ is a spinor field. The string equations, which 
are a generalization of the Einstein equations, can be written in the following form (see \cite{Stro}):
\bdm
\mbox{Ric}^g_{ij} - \frac{1}{4} H_{imn}H_{jmn} + 2 \cdot \nabla^g_i \partial_j
\Phi \ = \ 0, \quad \delta(e^{-2\Phi} H) \ = \ 0 \, .
\edm
The field equations are supplemented with the so-called Killing spinor equations
\bdm
(\nabla^g_X + \frac{1}{4} X \haken H) \cdot \Psi \ = \ 0, \quad (d\Phi - 
\frac{1}{2} H) \cdot \Psi \ = \ 0 \, .
\edm

Sometimes one requires that the $3$-form $H$ is closed, $dH = 0$, and 
solutions of this type are then called strong. The first of the Killing 
spinor equations suggests that the $3$-form $H$ should be the torsion form of 
a metric connection $\nabla$ with totally skew-symmetric torsion tensor $T=H$. 
Then the string equations and the Killing spinor equations can be written as follows (see \cite{IP}) :
\bdm
\mbox{Ric}^{\nabla} + \frac{1}{2}\delta(T) + 2 \cdot \nabla^gd\Phi=0, \quad 
\delta(T)=2 \cdot d\Phi^{\#} \haken T \, ,
\edm
\bdm
\nabla \Psi \ = \ 0, \quad 
(d\Phi - \frac{1}{2}T) \cdot \Psi \ = \ 0 \ .
\edm
$\mbox{Ric}^{\nabla}$ is the Ricci tensor of the metric connection $\nabla$, $\nabla^g$ is the Levi-Civita connection of the metric $g$ and $d\Phi^{\#}$ denotes the vector field dual to the 1-form $d\Phi$. If the dilation is constant, then the string equations are equivalent to the condition that the Ricci tensor 
of the connection vanishes \cite{IP} ,
\bdm
\mbox{Ric}^{\nabla}=0 \ , 
\edm
and the Killing spinor equations become
\bdm
\quad \nabla \Psi \ = \ 0, \quad T \cdot \Psi \ = \ 0 \ .
\edm

In particular, the spinor field is Riemannian harmonic. More general, the Riemannian Dirac operator $D^g$ acts on a $\nabla$-parallel spinor field
via the formula
\bdm
D^g \Psi + \frac{3}{4} \, T \cdot \Psi \ = \ 0 \ .
\edm

The number of preserved supersymmetries depends essentially on the number of $\nabla$-parallel spinors. In this paper we investigate solutions of all these field equations in the case of constant dilation in dimensions $n=5,6,7$. In dimension 7 any solution is strong (Theorem 5.4) and we derive a topological obstruction for the existence of such solutions (Remark 5.5).\\

A spinor parallel with respect to the Levi-Civita connection on a Riemannian manifold $M^n$ restricts its holonomy group (see \cite{H,Wa} and \cite{Mc}). 
In a similar way a $\nabla$-parallel spinor field reduces the structure group 
of the frame bundle. Conversely, let us start with a (non-integrable) $G$-structure on a Riemannian manifold and ask the question whether or not there exists a $G$-connection with a totally skew-symmetric torsion and at least one parallel spinor field. In dimension 3 the stabilizer of Spin(3) is trivial and therefore the connection is flat.  It is well known (see e.g. \cite{SSTV}) that in this case $(M^3,g,\nabla,T)$ carries (locally) a structure of a compact Lie group, $g$ is a biinvariant metric and $\nabla$ is the invariant connection with torsion given by the Lie bracket. In particular, on $SU(2)$ there exists at least one $\nabla$-parallel spinor.\\

The 4-dimensional case was investigated in earlier papers. The restricted
 holonomy group of $\nabla$ should be contained in SU(2) (see \cite{P1}) and this is equivalent to the local existence of a HKT structure, i.e., a hyperhermitian structure that is parallel with respect to $\nabla$ (see  \cite{IP}). Surprisingly, the geometry of $\nabla$ depends on the type of the parallel spinor (see \cite{DI}). If $M^4$ is compact, then the holonomy of $\nabla$ is contained in SU(2) if
and only if  $M^4$ is either a Calabi-Yau manifold or a Hopf surface (see \cite{IP}). We note that there exist Hopf surfaces that do not admit any (global) 
hyperhermitian structure although the holonomy of $\nabla$ is contained in 
SU(2) (see \cite{God}, \cite{IP}). These Hopf surfaces do not admit any 
$\nabla$-parallel spinors (see \cite{DI}), which shows that in the non-simply connected compact case the holonomy condition is not sufficient 
for the existence of $\nabla$-parallel spinors. \\

In higher dimensions we fix a subgroup $G$ of $Spin(n)$ preserving a spinor 
$\Psi_0$ as well as a geometric $G$-structure on a Riemannian manifold $(M^n,g)$. Then
\begin{itemize}
\item[a)] we describe the set of all $G$-connections with totally skew-symmetric torsion for any geometric type of $G$-structures. In particular, we decide whether or not a $G$-structure admits a connection with totally skew-symmetric torsion $T$ and derive a formula for the torsion;
\item[b)] We use the parallel spinor $\Psi_0$ and the algebraic properties of 
the Clifford multiplication in special dimensions in order to derive the field equation; 
\item[c)] We study the space of all $\nabla$-parallel spinors and compare it to the space of $\nabla$-harmonic spinors. In particular, we decide whether or not there exist $\nabla$-parallel spinor fields $\Psi$ such that $T \cdot \Psi 
= 0$.
\end{itemize}

In the second and third step we need the Schr\"odinger-Lichnerowicz-formula for a metric connection with totally skew-symmetric torsion. In particular, we
prove this formula in full generality, including the computation of the
curvature term. \\

We will give a complete answer in dimension $n=7$ and for  $G_2$-structures. In odd dimensions (almost contact metric structures) and in even dimensions 
(almost hermitian structures) we solve the first problem. However, these geometric structures reduce the structure group of the frame bundle only to the 
subgroup $U(k)$ , which does not coincide with the isotropy group of a spinor. Consequently, the connection $\nabla$ adapted to the geometric structure under consideration does not admit a $\nabla$-parallel spinor automatically and we obtain a further curvature condition for the existence of such spinors.  We 
investigate this condition and prove vanishing theorems for $\nabla$-harmonic spinors. 
%-----------------------------------------------------------------------------
\section{The curvature of connections with totally skew-symmetric torsion}
%-----------------------------------------------------------------------------
%
In this section we recall some notions concerning the curvature of a
metric connection with totally skew-symmetric torsion from \cite{IP}. Let $(M^n,g, \nabla, T)$ be an $n$-dimensional Riemannian manifold with
a metric connection $\nabla$ of totally skew-symmetric torsion $T$. The torsion tensor measures the difference between the connection $\nabla$ and the Levi-Civita connection $\nabla^g$ :
\bdm
g(\nabla_X Y,Z) \ = \ g(\nabla^g_X Y,Z) + \frac{1}{2} T(X,Y,Z) \ .
\edm
Let us fix some notation. The differential of an exterior form $\alpha$ 
is given by the formula
\bdm
d\alpha \ = \ \sum\limits^n_{i=1} e_i \wedge \nabla^g_{e_i}\alpha \, .
\edm
The codifferential of the form $\alpha$ can be calculated using either the 
Levi-Civita connection $\nabla^g$ or the connection $\nabla$ :
\bdm
\delta^g(\alpha) \ = \ - \,\sum \limits^n_{i=1} e_i \haken \nabla^g_{e_i}
\alpha, \quad \delta^{\nabla}(\alpha) \ =\ - \, \sum \limits^n_{i=1} e_i 
\haken \nabla_{e_i}\alpha \, .
\edm
On the $3$-form $T$, the two codifferentials coincide :
\bdm
\delta^g(T) \ = \ \delta^{\nabla}(T) \, .
\edm
This formula is a consequence of the assumption that $\nabla$ has a totally 
skew-symmetric torsion tensor. Let us introduce the $4$-form 
$\sigma^T$
\begin{eqnarray*}
\sigma^T (X,Y,Z,V) \hspace{-3mm} &:=& \hspace{-3mm} g(T(X,Y),T(Z,V))+
g(T(Y,Z), T(X,V))+ g (T(Z,X), T(Y,V)) \\
&=& \frac{1}{2} \sum\limits^n_{i=1} (e_i \haken T) \wedge (e_i \haken T)
(X,Y,Z,V) \, .
\end{eqnarray*}
Then the exterior derivative $dT$ of the torsion tensor $T$ is given in 
terms of $\nabla$ by the following formula (see e.g.\cite{IP})
\bdm
dT (X,Y,Z,V) \ = \ \sigma_{XYZ} \{ (\nabla_XT)(Y,Z,V)\} - (\nabla_V
T)(X,Y,Z)+ 2 \sigma^T(X,Y,Z,V) , 
\edm
where $\sigma_{XYZ}$ denotes the cyclic sum over $X,Y,Z$. Moreover, the curvature 
tensors of the Levi-Civita connection and the connection $\nabla$ are related 
via the formula 
\begin{eqnarray*}
R^g (X,Y,Z,V) &=& R^{\nabla}(X,Y,Z,V) - \frac{1}{2} (\nabla_X T)(Y,Z,V) +
\frac{1}{2} (\nabla_Y T)(X,Z,V)\\
&& - \frac{1}{4}g\Big(T(X,Y),T(Z,V)\Big) -\frac{1}{4}\sigma^T(X,Y,Z,V)\,.
\end{eqnarray*}
The first Bianchi identity for $\nabla$ can be written in the form
\bdm
\sigma_{XYZ} R^{\nabla}(X,Y,Z,V) \ = \ dT(X,Y,Z,V) - \sigma^T(X,Y,Z,V)+
(\nabla_V T)(X,Y,Z)  
\edm
and the difference of the Ricci tensors involves the codifferential of $T$ : 
\bdm
\mbox{Ric}^g(X,Y) \ = \ \mbox{Ric}^{\nabla}(X,Y) + \frac{1}{2} \delta^g
(T)(X,Y) - \frac{1}{4} \sum \limits^n_{i=1}g(T(e_i,X),T(Y,e_i)) \, .
\edm
In particular, the skew-symmetric part of the Ricci tensor 
$\mbox{Ric}^{\nabla}$ of $\nabla$ is given by the codifferential of the 
torsion tensor only :
\bdm
\mbox{Ric}^{\nabla}(X,Y) - \mbox{Ric}^{\nabla}(Y,X)\ = \ - \, \delta^g(T)(X,Y)\, .
\edm
We denote the scalar curvature of $\nabla$ by $\mbox{Scal}^{\nabla}$, i.e.,
\bdm
\mbox{Scal}^{\nabla} \ = \ \sum\limits^n_{i,j=1} R^{\nabla}
(e_i , e_j, e_j, e_i) \, .
\edm
%
%----------------------------------------------------------------------------
%----------------------------------------------------------------------------
\section{The Schr\"odinger-Lichnerowicz-formula for connections with totally 
skew-symmetic torsion}
%----------------------------------------------------------------------------
%
Consider an $n$-dimensional Riemannian spin manifold $(M^n,g, \nabla, T)$ 
with a metric connection $\nabla$ of totally skew-symmetric torsion $T$ and 
denote by $\Sigma M^n$ the spinor bundle. The Dirac operator $D$ depending on 
the connection $\nabla$ is defined by
\bdm
D \Psi \ = \ \sum\limits^n_{i=1} e_i \cdot \nabla_{e_i} \Psi , \quad
\Psi \in \Sigma M^n , 
\edm
where $e_1  \ldots e_n$ is an orthonormal basis. The Dirac operator $D$ is a 
formally selfadjoint operator since the torsion of the connection is totally 
skew-symmetric (see~\cite{FS}). In case of a Riemannian manifold the well-known Schr\"odinger-Lichnerowicz-formula expresses the square of the Dirac operator with respect to the Levi-Civita connection by the spinorial Laplace operator 
and some curvature term (see~\cite{Schroe,F1}). In the articles \cite{Bis,AT},
 a generalization of this formula for connections with arbitrary torsion is indicated. For connections with totally skew-symmetric torsion we shall derive the 
curvature term  and prove the following explicit formula : \\

{\bf Theorem 3.1. (S-L-formula)} {\it Let $(M^n,g, \nabla, T)$ be an
$n$-dimensional Riemannian spin manifold with a metric connection $\nabla$ of 
 totally skew-symmetric torsion $T$. Then, for any spinor field $\Psi$, 
the formula}
\bdm
D^2 \Psi = \nabla^* \nabla \Psi + \frac{3}{4} dT \cdot \Psi - \frac{1}{2} 
\sigma^T \cdot \Psi + \delta^g(T) \cdot \Psi - \sum\limits^n_{k=1} e_k 
\haken T \cdot \nabla_{e_k} \Psi + \frac{1}{4} \mbox{Scal}^{\nabla}\cdot 
\Psi 
\edm
{\it holds, where $\nabla^* \nabla$ is the Laplacian of $\nabla$ acting on 
spinors by}
\bdm
\nabla^* \nabla{\Psi} \ = \ - \sum\limits^n_{i=1} \nabla_{e_i} 
\nabla_{e_i} \Psi + \nabla_{\nabla_{e_i}^g{e_i}} \Psi\, .
\edm
{\bf Proof.} At a fixed point $p \in M^n$ we choose an orthonormal basis 
$e_1,\, \ldots, e_n$ such that  $(\nabla_{e_i} e_j)_p =0$ and 
$ [e_i, e_j]_p =-T(e_i, e_j)_p$. Then the vector field $e_i$ is 
$\nabla^g$-parallel in the direction of $e_i$ at the 
point $p$ : $(\nabla^g_{e_i} e_i)_p=0$. We calculate
\begin{eqnarray*}
D^2 \Psi \hspace{-2mm}&=& \hspace{-4mm}\sum\limits^n_{i,j=1} e_i \nabla_{e_i} e_j \nabla_{e_j} \Psi=- \sum\limits^n_{i=1} \nabla_{e_i} \nabla_{e_i} \Psi + \sum\limits^n_{i,j=1}R^{\nabla}(e_i, e_j) e_i \cdot e_j \cdot \Psi - \sum\limits^n_{i=1} e_i \haken T 
\cdot \nabla_{e_i} \Psi {} \\
&=& \nabla^* \nabla \Psi + \frac{1}{2}  \sum\limits_{i<j,k<l}
R^{\nabla}(e_i, e_j, e_k, e_l) e_i \cdot e_j \cdot e_k \cdot e_l \cdot \Psi - 
\sum\limits^n_{i=1} e_i \haken T \cdot \nabla_{e_i} \Psi \ .
\end{eqnarray*}
The curvature term in the latter equation can be written in the form
\bdm
\hat{\sigma}(R^{\nabla}) \cdot \Psi - \sum\limits_{j<k}
[\mbox{Ric}^{\nabla} (e_j, e_k) - \mbox{Ric}^{\nabla} (e_k, e_j)]e_j 
\cdot e_k \cdot \Psi + \frac{1}{2} \mbox{Scal}^{\nabla} \cdot \Psi \, , 
\edm
where $\hat{\sigma} (R^{\nabla})$ is the 4-form given by 
\bdm
\hat{\sigma} (R^{\nabla}) (X,Y,Z,V) \ = \ \sigma_{XYZ} 
\{ R^{\nabla}(X,Y,Z,V) - R^{\nabla}(V,X,Y,Z) \} \, .
\edm
Using the formulas for $dT$ as well as the formula comparing the curvature tensors $R^g$ and $R^{\nabla}$  we obtain
\bdm
\hat{\sigma} (R^{\nabla}) \ = \ \frac{3}{2} dT - \sigma^T \, . 
\edm
Inserting the latter formula as well as the formula for the skew-symmetric part
of the Ricci tensor into the expression for $D^2\Psi$ yields the desired 
formula. \hfill $\Box$\\

{\bf Corollary 3.2.} {\it Let $\Psi$ be a parallel spinor with respect to 
$\nabla$. Then the following formulas hold:}
\begin{eqnarray*}
\frac{3}{4} dT \cdot \Psi - \frac{1}{2} \sigma^T \cdot \Psi + \frac{1}{2}
\delta^g(T) \cdot \Psi + \frac{1}{4} \mbox{Scal}^{\nabla} \cdot \Psi  \ = \ 0 \, ,
\end{eqnarray*}
\begin{eqnarray*}
\Big(\frac{1}{2} X \haken dT + \nabla_XT \Big) \cdot \Psi - \mbox{Ric}^{\nabla} (X) \cdot \Psi \ = \ 0 \, . 
\end{eqnarray*}

{\bf Proof.} The first formula follows directly from Theorem 3.1. We prove the second one by  contracting the well-known formula 
\bdm
0 \ = \ \nabla \nabla \Psi \ = \ \sum \limits^n_{i,j=1} R^{\nabla}(e_i,e_j)
\cdot e_i \cdot e_j \cdot \Psi 
\edm
and using the formulas relating the symmetrization of the curvature tensor
$\mbox{R}^{\nabla}$ and the derivative of the torsion form $T$.\hfill$\Box$\\

The next formula compares the action of the Dirac operator with the action of the torsion form on spinors.\\

{\bf Theorem 3.3.} {\it Let $(M^n,g, \nabla,T)$ be an $n$-dimensional 
Riemannian spin manifold with a metric connection $\nabla$ of totally 
skew-symmetric torsion. Then}
\bdm
DT+TD\ =\ dT + \delta^g(T) - 2 \cdot \sigma^T - 2 \sum\limits^n_{i=1} 
e_i \haken T  \cdot\nabla_{e_i} \, . 
\edm

{\bf Proof.} The proof is similar to the proof of Theorem 3.1. \hfill $\Box$\\

Let us apply Theorem 3.1 and Theorem 3.3 in case of a $\nabla$-harmonic spinor
field $\Psi$ on a compact manifold $M^n$. Since $\delta^g(T)$ is a $2$-form,
 the real part of the hermitian product $(\delta^g(T) \cdot \Psi,\Psi)$ vanishes. Moreover, the Dirac operator $D$ is symmetric. Combining Theorem 3.1 and Theorem 3.3 we obtain the condition
\bdm
\int_{M^n} \Big(\| \nabla \Psi \|^2 + \frac{1}{4}(dT \cdot \Psi, \Psi) + 
\frac{1}{2}(\sigma^T \cdot \Psi, \Psi) + \frac{1}{4} \mbox{Scal}^{\nabla}
\cdot \|\Psi \|^2 \Big) \ = \ 0 \ .
\edm
This formula proves the following vanishing theorem.\\

{\bf Theorem 3.4.} {\it Let $(M^n,g, \nabla, T)$ be a compact
Riemannian spin manifold with a metric connection $\nabla$ of
 totally skew-symmetric torsion $T$. Suppose, moreover, that the eigenvalues
of the endomorphism} $dT + 2 \cdot \sigma^T + \mbox{Scal}^{\nabla}$ {\it 
acting on spinors are non-negative. Then any $\nabla$-harmonic spinor field is 
$\nabla$-parallel. In case the eigenvalues of the endomorphism are positive,
 there are no $\nabla$-parallel spinor fields.}\\

%----------------------------------------------------------------------------
\section{$G_2$-connections with totally skew-symmetric torsion}
%----------------------------------------------------------------------------
%
We study under which conditions a fixed $G$-structure on a Riemannian manifold
admits an affine connection preserving the $G$-structure and having
totally skew-symmetric torsion tensor. For this purpose we describe the different geometric types of $G$-structures from the point of view of gauge theory using a certain $1$-form $\Gamma$ with values in the associated bundle of typical fibres $\son/  \altg$. This approach is completely equivalent to the classification of different geometric $G$-structures used in differential geometry and 
mainly based on the decomposition of the covariant derivative of the tensor related with the $G$-structure. The advantage of our approach is that the method 
applies even in cases where the $G$-structure is not defined by a tensor (see~\cite{tf-weak, Swann}). To begin with, let $(M^n,g)$ be an oriented Riemannian manifold and denote by ${\cal F} (M^n)$ its frame bundle. The Levi-Civita connection is a $1$-form
\bdm
 Z^g:T\big({\cal F}(M^n)\big) \longrightarrow \son
\edm
with values in the Lie algebra $\son$. Its torsion tensor vanishes. We fix a
closed subgroup $G$ of the orthogonal group $SO(n)$. A $G$-structure on 
$M^n$ is a $G$-subbundle ${\cal R} \subset {\cal F} (M^n)$. We decompose the 
Lie algebra $\son$ into the subalgebra $\altg$ and its orthogonal  complement 
$\altm$ :
\bdm
 \son \ = \ \altg \oplus \altm \ . 
\edm
In a similar way we decompose the restriction of the 1-form $Z^g$ 
\bdm
Z^g|_{T({\cal R})} \ = \ \tilde{Z} \oplus \Gamma \, . 
\edm
$\tilde{Z}$ is a connection in the principal $G$-bundle ${\cal R}$ and $\Gamma$
is a $1$-form with values in the associated bundle ${\cal R} \times_G \altm$. 
The different geometric types of $G$-structures are defined by the irreducible 
$G$-components of the representation ${\Bbb R}^n \otimes \altm$. An arbitrary 
$G$-connection $Z$ differs from $\tilde{Z}$ by an $1$-form
$\Sigma$ with values in the Lie algebra $\altg$,
\bdm
 Z \ = \ Z^g - \Gamma + \Sigma \, ,
\edm
and the corresponding covariant derivative $\nabla$ is given by the formula
\bdm
 \nabla_X Y \ = \ \nabla^g_X Y -  \Gamma (X)(Y) + \Sigma (X)(Y) \, . 
\edm
Since the Levi-Civita connection $\nabla^g$ is torsion free, the torsion tensor
of $\nabla$ depends on $\Gamma$ and $\Sigma$ :
\bdm
T(X,Y,Z)\, =\, - g\big(\Gamma (X)(Y),Z\big)+g\big(\Gamma (Y)(X),Z\big) 
 + g \big(\Sigma (X)(Y),Z\big) - g\big(\Sigma (Y)(X),Z\big) \, . 
\edm
$T$ is a $3$-form if and only if
\bdm
g\big(\Gamma(Y)(X),Z\big) + g\big(\Gamma (Z)(X),Y\big) \ = \ 
g \big(\Sigma (Z)(X), Y \big)+g \big(\Sigma (Y)(X),Z \big)
\edm
holds. Now we introduce the following $G$-invariant maps:
\begin{eqnarray*}
\Phi :& {\Bbb R}^n \otimes \altg \to {\Bbb R}^n \otimes S^2 
({\Bbb R}^n),& \ \Phi (\Sigma)(X,Y,Z) := g \big(\Sigma (Z)(X),Y \big)+
g \big(\Sigma (Y)(X),Z \big),\\
\Psi :& {\Bbb R}^n \otimes \altm \to {\Bbb R}^n \otimes S^2 ({\Bbb R}^n),& \
\Psi (\Gamma)(X,Y,Z) := g \big(\Gamma (Y)(X),Z \big)+
g \big(\Gamma (Z)(X),Y \big) \, .
\end{eqnarray*}
Consequently, we proved the following\\

{\bf Proposition 4.1.} {\it A $G$-reduction ${\cal R} \subset {\cal F} (M^n)$
admits a $G$-connection $Z$ with a totally skew-symme\-tric torsion tensor $T$ if and only if $\Psi (\Gamma)$ is contained in the image of the homomorphism 
$\Phi$. In this case the set of all these connections $Z$ is an affine space
over the vector space $\ker (\Phi)$.}\\

We will use representation theory in order to study the diagramme
$$
\begin{diagram}
\node{{\Bbb R}^n \otimes \altg} \arrow{e,t}{\Phi} \node{{\Bbb R}^n 
\otimes S^2 ({\Bbb R}^n)}\\
\node{{\Bbb R}^n \otimes \altm} \arrow{ne,b}{\Psi} \node{} 
\end{diagram}
$$
By splitting the $G$-representation ${\Bbb R}^n \otimes \altm$ into
irreducible components we can decide whether or not the image of a certain 
component is contained in the image of $\Phi$. In this way we characterize 
the geometric $G$-structures admitting a $G$-connection $Z$ with a
totally skew-symmetric torsion tensor. We will apply this general method to 
the subgroup $G_2 \subset SO(7)$. Therefore, let us recall some notions of $G_2$-geometry in dimension seven. The group 
$G_2$ is the isotropy group of the $3$-form in seven variables
\bdm
\omega^3 \ := \ e_1 \wedge e_2 \wedge e_7 +  e_1 \wedge e_3 \wedge e_5 - 
 e_1 \wedge e_4 \wedge e_6 -  e_2 \wedge e_3 \wedge e_6 -  e_2 \wedge e_4 \wedge e_5 +  e_3 \wedge e_4 \wedge e_7 +  e_5 \wedge e_6 \wedge e_7 \ .
\edm
The $3$-form $\omega^3$ corresponds to a real spinor $\Psi_0 \in \Delta_7$ and, therefore, $G_2$ can be defined as the isotropy group of a non-trivial real spinor. We identify the Lie algebra of the group $SO(7)$ with the space of all $2$-forms :
\bdm
\mbox{\graf so}(7) \ = \ \Lambda^2({\Bbb R}^7) \ = \ \Big\{ \sum_{i < j} \omega_{ij} \cdot e_i \wedge e_j \Big\} \ .
\edm 
The Lie algebra $\altg_2$ of the group $G_2$ is given by the equations
\bdm
\omega_{12}+\omega_{34}+\omega_{56} \ = \ 0, \quad - \, \omega_{13}+\omega_{24}-\omega_{67} \ = \ 0, \quad \omega_{14}+\omega_{23}+\omega_{57} \ = \ 0,
\edm
\bdm
\omega_{16}+\omega_{25}-\omega_{37} \ = \ 0, \quad \quad \omega_{15}-\omega_{26}-\omega_{47} \ = \ 0, \quad \omega_{17}+\omega_{36}+\omega_{45} \ = \ 0,
\edm
\bdm
\quad \omega_{27}+\omega_{35}-\omega_{46} \ = \ 0 \ .
\edm
The space ${\Bbb R}^7 := \Lambda^1_7$ is an irreducible $G_2$-representation, the $2$-forms $\Lambda^2 = \Lambda^2_7 \oplus \Lambda^2_{14}$ split into two irreducible $G_2$-components :
\begin{eqnarray*}
\Lambda^2_7&:=&\{\alpha^2 \in \Lambda^2 \ : \ *(\omega^3 \wedge \alpha^2) =
2 \cdot \alpha^2 \} \ = \ \{ X \haken \omega^3 \ : \ X \in {\Bbb R}^7 \} \ ,\\
\Lambda^2_{14}&:=&\{\alpha^2 \in \Lambda^2 \ : \ *(\omega^3 \wedge \alpha^2) = - \, \alpha^2 \} \ = \ \altg_2 \ .
\end{eqnarray*}
The space of $3$-forms $\Lambda^3 = \Lambda^3_1 \oplus \Lambda^3_7 \oplus \Lambda^3_{27} $ decomposes into three irreducible $G_2$-components :
\begin{eqnarray*}
\Lambda^3_1 &:=& \{t \cdot \omega^3 \ : \ t \in {\Bbb R} \} \ ,\\ 
\Lambda^3_7 &:=& \{ *(\omega^3 \wedge \alpha^1) \ : \ \alpha^1 \in \Lambda^1\} \ = \ \{ X \haken *\omega^3 \ : \ X \in {\Bbb R}^7 \} \ , \\
\Lambda^3_{27} &:=& \{ \alpha^3 \in \Lambda^3 \ : \ \alpha^3 \wedge \omega^3
=0, \ \alpha^3 \wedge *\omega^3=0 \} \ .
\end{eqnarray*}
The representation $\Lambda^3_{27}$ is isomorphic to the representation of
$G_2$ in the space $S^2_0({\Bbb R}^7)$ of all traceless symmetric bilinear forms. \\

{\bf Proposition 4.2.} {\it The map \
$\Phi : {\Bbb R}^7 \otimes \altg_2 \rightarrow {\Bbb R}^7 \otimes S^2
({\Bbb R}^7)$ \ is injective.}\\

{\bf Corollary 4.3.} {\it Let $(M^7, g, \omega^3)$ be an oriented, 
$7$-dimensional Riemannian manifold with a fixed $G_2$-structure $\omega^3$. 
Then there exists at most one affine connection $\nabla$ such that $\nabla
\omega^3 =0$ and the torsion tensor $T$ is a $3$-form.}\\

{\bf Proof.} Given $\Sigma \in {\Bbb R}^7 \otimes \altg_2$, the condition
$\Phi (\Sigma)=0$ is equivalent to
\bdm
 Z \haken \Sigma (Y) + Y \haken \Sigma (Z) \ = \ 0
\edm
for any two vectors $Y,Z \in {\Bbb R}^7$. Using the standard basis 
$e_{\alpha} \wedge e_{\beta}$ of the Lie algebra $\mbox{\graf so}(7)$
we decompose the elements $\Sigma (e_i) \in \altg_2 \subset \mbox{\graf so}(7)$,
\bdm
 \Sigma (e_i) \ :=\ \sum\limits_{1 \le \alpha, \beta \le 7} \omega_{i
\alpha \beta} \cdot e_{\alpha} \wedge e_{\beta} \, . 
\edm
The condition $\Phi (\Sigma)=0$ implies $\omega_{i \alpha \beta} =
- \omega_{\alpha i \beta}$. Consequently, $\Sigma$ depends on 35 parameters 
$\omega_{i \alpha \beta} \ (1 \le i < \alpha < \beta \le 7)$. Moreover, $\Sigma (e_i)$ belongs to the Lie-Algebra $\altg_2$ and hence we obtain seven 
equations for any index $i=1,2, \ldots , 7$. Altogether, these are 49 equations for 35 variables and a careful examination of the system yields the result that $\omega_{i \alpha 
\beta} =0$ is the only solution. \mbox{} \hfill $\Box$\\

The low-dimensional $G_2$-representations, their highest weights, dimensions
etc. are listed in the following table (see~\cite{FH}):\\

\begin{center}
\begin{tabular}{|c|c|c|}
\hline &&\\
highest weight & dimension & space  \\ && \\ \hline
&&\\
$(0,0)$ & $1$ & $\Lambda^0_1 := \Lambda^0 ({\Bbb R}^7)= \Lambda^3_1$\\
&&\\
$(1,0)$ & $7$ & \mbox{}\quad  \quad $\Lambda^1_7 := \Lambda^1 ({\Bbb R}^7)= 
\Lambda^2_7= \Lambda^3_7$ \quad \quad  \mbox{}\\
&&\\
$(0,1)$ & $14$ & $\altg_2 = \Lambda^2_{14}$\\
&&\\
$(2,0)$ & $27$ & $\Lambda^3_{27} = S^2_0 ({\Bbb R}^7)$\\
&&\\
$(1,1)$ & $64$ & $\Lambda_{64}$\\
&&\\
$(3,0)$ & $77$ & $\Lambda_{77}$\\ &&\\ \hline
\end{tabular}
\end{center}

We now compute the $G_2$-decomposition of the three representations related 
with the maps $\Psi$ and $\Phi$.\\

{\bf Proposition 4.4.} {\it The following decompositions into irreducible
$G_2$-representations hold:}
\begin{enumerate}
\item[1)] ${\Bbb R}^7 \otimes \altm = \Lambda^0_1 \oplus \Lambda^1_7 \oplus 
\Lambda^2_{14} \oplus \Lambda^3_{27}$;
\item[2)] ${\Bbb R}^7 \otimes \altg_2 = \Lambda^1_7 \oplus \Lambda^3_{27} 
\oplus \Lambda_{64}$;
\item[3)] ${\Bbb R}^7 \otimes S^2 ({\Bbb R}^7)= 2 \Lambda^1_7 \oplus \Lambda^2_{14} \oplus \Lambda^3_{27} \oplus \Lambda_{64} \oplus \Lambda_{77}$.
\end{enumerate}

{\bf Proof:} The first decomposition is an elementary one, the second and 
third decomposition can be obtained using a suitable computer programme. \mbox{}
\hfill $\Box$\\

\newcommand{\rs}{{\Bbb R}^7}
\newcommand{\rn}{{\Bbb R}^n}

Since the $G_2$-map $\Phi : \rs \otimes \altg_2 \to \rs \otimes S^2 (\rs)$ is
injective and the multiplicity of $\Lambda^3_{27}$ in $\rs \otimes S^2 (\rs)$
is one, we obtain\\

{\bf Corollary 4.5.}
\begin{enumerate}
\item[1)] $\Psi (\Lambda^0_1 \oplus \Lambda^3_{27}) \subset \mbox{Im} (\Phi) $;
\item[2)] $\Psi (\Lambda^2_{14}) \cap \mbox{Im} (\Phi) =0$. 
\end{enumerate}

Finally, we have to decide whether or not the space $\Psi (\Lambda^1_7)$
is contained in the image of $\Phi$. Since the representation $\Lambda^1_7$ 
has multiplicity two in $\rs \otimes S^2 (\rs)$, we cannot use a universal
argument is before.\\

{\bf Proposition 4.6.} {\it $\Psi (\Lambda^1_7)$ is contained in}
$\mbox{Im} (\Phi)$.\\

{\bf Proof.} First of all we compute $\Psi (\Gamma)$ for a given vector
$\Gamma \in \Lambda^1_7$. The element in $\rs \otimes \rs$ related to
$\Gamma$ is the 2-form (skew-symmetric endomorphism) $\Gamma \haken \omega^3$
and therefore we obtain
\
\bdm
 \Psi (\Gamma)(X,Y,Y) \ = \ 2 \sum\limits^7_{i=1} \omega^3 (\Gamma , Y,
e_i) \cdot \omega^3 (e_i, X,Y) \ = \ 2 g(\Gamma,X) \cdot g(Y,Y) - 2 g(\Gamma,X)\cdot g(X,Y) \, . 
\edm
Suppose now that $\Psi (\Lambda^1_7) \subset \mbox{Im}(\Phi)$. Since
$\Phi$ is injective, there should exist a map $\Sigma : \Lambda^1_7 \to \rs
\otimes  \altg_2$ such that
\bdm
\Phi \big(\Sigma(\Gamma) \big) \ = \ \Psi (\Gamma) \, .
\edm
The multiplicity of $\Lambda^1_7$ in $\rs \otimes \altg_2$ equals one.
Consequently, $\Sigma (\Gamma)$ is proportional to the map
\bdm
 \Sigma_0 (\Gamma)(Y) \ := \ \mathrm{pr}_{\altg_2} (\Gamma \wedge Y) \, ,
\edm
where $\mathrm{pr}_{\altg_2}: \Lambda^2 (\rs) = \mbox{\graf so}(7) \to  
\altg_2$ is the orthogonal projection. The $2$-forms 
\bdm
\frac{1}{\sqrt{3}} e_1 \haken \omega^3 ,\quad \ldots \quad ,
\frac{1}{\sqrt{3}} e_7 \haken \omega^3
\edm
constitute an orthonormal basis of the subspace 
$\altm \subset \mbox{\graf so}(7)$ and the projection is given by the formula
\bdm
\mbox{pr}_{\altg_2}(\alpha^2)\ = \ \alpha^2 - \frac{1}{3} \sum\limits^7_{i=1} 
( e_i \haken \omega^3 , \alpha^2) \cdot (e_i \haken \omega^3) \, . 
\edm
Now we compute $\Phi \big(\Sigma_0 (\Gamma)\big)$:
\begin{eqnarray*}
\Phi \big(\Sigma_0 (\Gamma)\big)(X,Y,Y) &=& 2 \cdot g(X \haken 
\mathrm{pr}_{\altg_2} (\Gamma \wedge Y),Y)\\
&=& 2 \cdot (\Gamma \wedge Y)(X,Y) - \frac{2}{3} \sum\limits^7_{i=1}
(e_i \haken \omega^3 , \Gamma \wedge Y) \cdot \omega^3 (e_i, X, Y) \, . 
\end{eqnarray*}
Finally, we obtain
\bdm
 \Phi \big(\Sigma_0 (\Gamma)\big)(X,Y,Y) \ = \ \frac{2}{3} \Psi (\Gamma)(X,Y,Y) \, , 
\edm
and the latter formula proves the Proposition. \mbox{} \hfill $\Box$\\

To summarize the previous discussion we proved the following theorem.\\

{\bf Theorem 4.7.} {\it Let $(M^7,g, \omega^3)$ be a $7$-dimensional
Riemannian manifold with a $G_2$-structure $\omega^3$. The following
conditions are equivalent:
\begin{itemize} 
\item[$1)$] 
The $\Lambda^2_{14}$-component of $\Gamma$ is zero.
\item[$2)$]
There exists an affine connection $\nabla$ with totally skew-symmetric torsion
preserving the $G_2$-structure,
\bdm
\nabla \omega^3 \ = \ 0 \, .
\edm
\end{itemize}
In this case the connection $\nabla$ is unique.}\\

The Riemannian covariant derivative $\nabla^g \omega^3$ of the $G_2$-structure 
is given by the formula
\bdm
 \nabla^g \omega^3 \ = \ \varrho_3 (\Gamma)(\omega^3) \, ,
\edm
where $\varrho_3$ is the differential of the representation of $SO(7)$ in 
$\Lambda^3 (\rs)$. The $G_2$-map of $\rs \otimes \altm$ into $\rs \otimes
\Lambda^3 (\rs)$ given by the formula $X \otimes Y \to X \otimes \varrho_3
(Y) \omega^3$ is injective (see~\cite{tf-weak}). Consequently, the different 
geometric classes of $G_2$-structures introduced by Fernandez/Gray (see~\cite
{FG}) using the covariant derivative $\nabla^g \omega^3$ can be defined via the
algebraic type of $\Gamma \in \Lambda^1 (M^7) \otimes \Lambda^1 (M^7)$. The
condition $\Gamma \in \Lambda^0_1 \oplus \Lambda^1_{7} \oplus \Lambda^3_{27}$ 
characterizes the so-called integrable $G_2$-structures, 
i.e., the $G_2$-structures of type $W_1 \oplus W_4 \oplus W_3$ in the notation 
of \cite{FG}. Let us decompose $\Gamma$ into its three parts:
\bdm
\Gamma \ = \ \lambda \cdot \mbox{Id}_{T(M^7)} \oplus \beta \oplus \Gamma^3_{27} \, .
\edm
The $\Lambda^0_1$ -part acts as a map
$\lambda \cdot \mbox{Id}_{T(M^7)}: \rs \to \altm$ via the formula
\bdm
(\lambda \cdot \mbox{Id}_{T(M^7)})(X) \ := \ \frac{\lambda}{12} \cdot 
(X \haken \omega^3) \, .
\edm
The component $\beta \in \Lambda^1_7$ is a vector field and we will use the
embedding $\rs \subset \rs \otimes \altm$ given by the equation
\bdm
\beta(X) \ := \ \frac{1}{4} \cdot \mbox{pr}_{\altm}(\beta \wedge X) \ = \ 
\frac{1}{12}\sum^7_{i=1} (\beta \wedge X, e_i \haken \omega^3) \cdot (e_i \haken \omega^3) \, .
\edm
The third component $\Gamma^3_{27} \in \Lambda^3_{27}$ is a $3$-form defining
a map $\Gamma^3_{27}: \rs \to \altm$ by the formula
\bdm
\Gamma^3_{27}(X) \ := \ \frac{1}{2} \cdot \mbox{pr}_{\altm}(X \haken \Gamma^3_{27})
\ = \ \frac{1}{6}\sum^7_{i=1} (X \haken \Gamma^3_{27}, e_i \haken \omega^3) 
\cdot (e_i \haken \omega^3) \, .
\edm
We describe the action $\varrho_3(\alpha^2)$ of a 
$2$-form $\alpha^2 \in \Lambda^2 (\rs)= \mbox{\graf so} (7)$ on the form
$\omega^3$. Suppose that the projection of $\alpha^2$ onto the space
$\Lambda^2_7$ is given by a vector $Z \in \rs$, 
\bdm
\mbox{pr}_{\mbox{\graf m}}(\alpha^2) \ = \ Z \haken \omega^3 \, .
\edm
Then the representation $\varrho_3 (\alpha^2)$ acts on the 3-form $\omega^3$ 
by the formula (see~\cite{FKMS}) 
\bdm
\varrho_3 (\alpha^2)(\omega^3) \ = \ - \ 3 \cdot (Z \haken * \omega^3) \, .
\edm
Consequently, we obtain a formula for the covariant derivative 
$\nabla^g \omega^3$ involving the function $\lambda$, the vector field 
$\beta$ as well as the $3$-form $\Gamma^3_{27}$:
\bdm
\nabla^g_X \omega^3  \ = \ - \, \frac{\lambda}{4} \cdot (X \haken * \omega^3) 
- \sum^7_{i=1} \big(\frac{1}{4} \beta \wedge X + \frac{1}{2} X \haken 
\Gamma^3_{27} , e_i \haken \omega^3) \cdot (e_i \haken * \omega^3) \, .
\edm
Now we can express the differential $d \omega^3$ and the codifferential 
$\delta^g(\omega^3)$ of the $3$-form by $\Gamma$. For example, we have
\bdm
\delta^g(\omega^3) \ = \ - \, \sum^7_{j=1} e_j \haken \nabla^g_{e_j}
 \omega^3 \ = \  \sum^7_{i,j=1} \Big( \frac{1}{4}\beta \wedge e_j + \frac{1}{2} e_j \haken 
\Gamma^3_{27}, e_i \haken \omega^3 \Big) \cdot \Big(e_j \haken (e_i \haken * 
\omega^3)\Big)\,. 
\edm
The map 
\bdm
\Gamma^3_{27} \in \Lambda^3_{27} \longmapsto \sum^7_{i,j=1} (e_j \haken 
\Gamma^3_{27}, e_i \haken \omega^3) \cdot (e_j \haken (e_i \haken * \omega^3)) \in \Lambda^2(\rs) = \Lambda^2_7 \oplus \Lambda ^2_{14}
\edm
is obviously trivial and the map 
\bdm
\beta \in \rs \longmapsto \sum^7_{i,j=1} (\beta \wedge e_j, e_i \haken 
\omega^3) \cdot (e_j \haken (e_i \haken * \omega^3)) \in \Lambda^2(\rs) = 
\Lambda^2_7 \oplus \Lambda ^2_{14}
\edm 
must be proportional to the map $ \beta \to \beta \haken \omega^3$ since the 
multiplicity of the $7$-dimensional component in $\Lambda^2(\rs)$ is one. 
Computing the constant we obtain
\bdm
\sum^7_{i,j=1} (\beta \wedge e_j, e_i \haken \omega^3) \cdot (e_j \haken (e_i \haken * \omega^3)) \ = \ - \, 4 \cdot (\beta \haken \omega^3)
\edm
and, finally,
\bdm
\delta^g(\omega^3) \ = \ - \, \beta \haken \omega^3 \, .
\edm
We handle the differential $d \omega^3$ in a similar way. Using again that 
the multiplicities of $\rs$ and $\Lambda^3_{27}$ in the $G_2$-representation
$\Lambda^4(\rs)$ are one we verify the idendities
\bdm
\sum^7_{i,j=1} (e_j \haken \Gamma^3_{27}, e_i \haken \omega^3) \cdot (e_j 
\wedge (e_i \haken * \omega^3)) \ = \ - \, 2 \cdot (* \Gamma^3_{27})\,,
\edm
\bdm
\sum^7_{i,j=1} (\beta \wedge e_j , e_i \haken \omega^3) \cdot (e_j \wedge (e_i \haken * \omega^3)) \ = \ - \, 3 \cdot (\beta \wedge \omega^3) \, .
\edm
Then we obtain the formula
\bdm
d \omega^3 \ = \ \sum^7_{i=1} e_i \wedge \nabla^g_{e_i} \omega^3 \ = \ 
- \lambda \cdot (*\omega^3) + * \Gamma^3_{27} + \frac{3}{4} \cdot (\beta 
\wedge \omega^3) \, .
\edm 
We compute now the $1$-form $\Sigma(\Gamma) \in \rs \otimes \altg_2$ defined by the condition $\Phi(\Sigma(\Gamma)) = \Psi(\Gamma)$. Since $\Psi(\lambda 
\cdot \mbox{Id}) = 0$ we have  $\Sigma(\lambda \cdot \mbox{Id}) = 0$, i.e., $\Sigma(\Gamma)$ does not depend on the $\Lambda^0_1$ -part of $\Gamma$. The 
$1$-form $\Psi ({\Gamma})$ with values in $S^2({\Bbb R}^7)$ is given by
\bdm
\Psi ({\Gamma})(X,Y,Y) \ = \ 2 \cdot g \Big({\Gamma} (Y)(X),Y\Big) \ = \ 
g\Big(X \haken \mbox{pr}_{\altm}\big(Y \haken \Gamma^3_{27} + \frac{1}{2} \beta \wedge Y \big),Y\Big) \, . 
\edm
On the other hand, let us introduce the map $\Sigma: \Lambda^1_7 \oplus 
\Lambda^3_{27} \to \rs \otimes \altg_2$ given by the formula
\bdm
\Sigma (\Gamma)(Y) \ := \ - \, \frac{1}{2}\cdot \mbox{pr}_{\altg_2} \Big(Y \haken 
\Gamma^3_{27} - \frac{1}{4} \beta \wedge Y \Big) \, .
\edm
Proposition 4.6 and a direct calculation yields that $\Sigma(\Gamma)$ is 
indeed the $1$-form with values in the Lie algebra $\altg_2$ corresponding 
to $\Gamma$:
\begin{eqnarray*}
\Phi \big(\Sigma (\Gamma^3_{27})\big)(X,Y,Y) &=& - \, g \big(X \haken \mbox{pr}_{\altg_2}(Y \haken \Gamma^3_{27}), Y\big)  \\
&=& - \, (Y \haken \Gamma^3_{27})(X,Y) + g \big(X \haken \mbox{pr}_{\altm}(Y \haken \Gamma^3_{27}), Y\big) \\
&=& \Psi (\Gamma^3_{27})(X,Y,Y) \, .
\end{eqnarray*}
We introduce a new $1$-form $\Gamma^* (X):= \Gamma(X) - \Sigma(\Gamma(X))$,
\bdm
\Gamma^*(X) \ = \ \frac{\lambda}{12} \cdot (X \haken \omega^3) + \frac{1}{2}
\cdot (X \haken \Gamma^3_{27}) + \frac{3}{8} \cdot \mbox{pr}_{\altm}(\beta \wedge X) - \frac{1}{8}\cdot (\beta \wedge X) 
\edm
and remark that the torsion form $T$ is given by 
\bdm
T(X,Y,Z) \ = \ - \, g(\Gamma^*(X)(Y),Z) +  g(\Gamma^*(Y)(X),Z) \, .
\edm
Consequently, we obtain a formula for the torsion form $T$,
\bdm
T \ = \ - \, \frac{\lambda}{6} \cdot \omega^3 - \Gamma^3_{27} + T_{\beta} \, ,
\edm
where the $3$-form $T_{\beta}$ is defined by the equation
\begin{eqnarray*}
T_{\beta}(X,Y,Z) & := & \frac{3}{8} \cdot \Big(\mbox{pr}_{\altm}(\beta \wedge Y)(X,Z) - \mbox{pr}_{\altm}(\beta \wedge X)(Y,Z)\Big)\\
& + & \frac{1}{8}\cdot \Big( g(\beta,Y)\cdot g(X,Z) - g(\beta,X) \cdot g(Y,Z)
\Big) \, .
\end{eqnarray*}
The map $\beta \in \rs \longmapsto T_{\beta} \in \Lambda^3(\rs)$ must be a multiple of the map $\beta \longmapsto \beta \haken * \omega^3$ since the multiplicity of $\rs$ in $\Lambda^3(\rs)$ is one. Computing algebraically the constant 
we obtain the equation
\bdm
T_{\beta} \ = \ - \, \frac{1}{4}\cdot ( \beta \haken * \omega^3) \, .
\edm
We thus computed the torsion form of the unique connection preserving the
$G_2$-structure:
\bdm
T \ = \ - \, \frac{\lambda}{6} \cdot \omega^3 - \Gamma^3_{27} - \frac{1}{4}
\cdot ( \beta \haken * \omega^3) \, .
\edm
Using the equations for the exterior differential and the codifferential of the form $\omega^3$ we can substitute the function $\lambda$ as well as the $3$-form $\Gamma^3_{27}$:
\bdm
\lambda \ = \ -\,\frac{1}{7}\cdot (d \omega^3,*\omega^3), \quad \Gamma^3_{27} \ = \ * d \omega^3 + \lambda \cdot \omega^3 - \frac{3}{4} *(\beta \wedge \omega^3) \, .
\edm
For any vector $\beta$ we have
\bdm
*(\beta \wedge \omega^3) \ = \ - \, (\beta \haken *\omega^3)\, .
\edm
Let us summarize the result.\\

{\bf Theorem 4.8.} {\it Let $(M^7,g, \omega^3)$ be a $7$-dimensional 
$G_2$-manifold of type $\Lambda^0_1 \oplus \Lambda^1_7 \oplus 
\Lambda^3_{27}$. The torsion form of the unique affine connection $\nabla$ 
preserving the structure with totally skew-symmetric torsion is given by the formula}
\bdm
T  \ = \ \frac{1}{6}\cdot (d \omega^3, * \omega^3) \cdot \omega^3 - \, * d \omega^3
+  *(\beta \wedge \omega^3) \, .
\edm
{\it The vector field $\beta$ as well as the differential and the codifferential of the $3$-form $\omega^3$ are related by}
\bdm
\delta^g(\omega^3) \ = \ - \, (\beta \haken \omega^3), \quad 
d\omega^3 \ = \ \frac{1}{7}\cdot(d\omega^3, * \omega^3) \cdot (* \omega^3) + * \Gamma^3_{27} + \frac{3}{4} \cdot (\beta \wedge \omega^3) \ .
\edm
\\
The particular case $\Gamma = \lambda \cdot \mbox{Id}_{T(M^7)}$ 
corresponds to nearly parallel $G_2$-structures. In this case, the Riemannian 
manifold $(M^7,g)$ is Einstein, the parameter $\lambda$ is constant and related to the scalar curvature of $M^7$ (see~\cite{FKMS}).\\

{\bf Corollary 4.9.} {\it Let $(M^7,g, \omega^3)$ be a $7$-dimensional 
Riemannian manifold with a nearly parallel $G_2$-structure} 
$(\Gamma = \lambda \cdot \mbox{Id})$. {\it Then there exists a unique affine 
connection $\nabla$ such that}
$$\nabla \omega^3 \ = \ 0 \quad \mathit{and} \quad  T \ \mathit{is\ a\ }  
3 \mbox{\it -form.}$$
{\it The torsion tensor is given by the formula} $6 \cdot T = (d \omega^3,
*\omega^3) \cdot \omega^3$. $T$ {\it is} $\nabla$-{\it parallel and coclosed,} $\nabla T = \delta T = 0$. \\

{\bf Corollary 4.10.} {\it Let $(M^7,g,\omega^3)$ be a $7$-dimensional nearly
parallel $G_2$-manifold. Then the triple $(M^7,g,T^*:= 3 \cdot T)$ is a solution of the string equations with constant dilation:}
\bdm
\mbox{Ric}^g_{ij} - \frac{1}{4} T^*_{imn}T^*_{jmn} \ = \ 0 , \quad \delta^g(T^*) \ = \ 0 \,.
\edm
A {\it cocalibrated} $G_2$-structure is defined by the condition that $\omega^3$ is coclosed, $\delta^g(\omega^3) = 0$. Equivalently, $\Gamma$ depends only on a 
function $\lambda$ and on a $3$-form $\Gamma^3_{27}$ of type $\Lambda^3_{27}$,
\bdm
\Gamma \ = \ \lambda \cdot \mbox{Id} + \Gamma^3_{27} \, .
\edm
The differential as well as the torsion form $T$ 
of cocalibrated  $G_2$-structures are given by the formulas:
\bdm
d \omega^3 \ = \ - \,  \lambda \cdot (* \omega^3) + (* \Gamma^3_{27}), \quad 
\lambda \ = \ - \, \frac{1}{7}\cdot (d \omega^3, * \omega^3) \, , 
\edm
\bdm
T \ = \ - \, (*d \omega^3) - \frac{7}{6} \cdot \lambda \cdot \omega^3, \quad 
d*T \ = \ - \, \frac{7}{6} \cdot d\lambda \wedge * \omega^3 \ .
\edm
In particular, $T$ is coclosed if and only if $\lambda$ is constant. This occurs, for example, if the $G_2$-structure is of pure type $\Lambda^0_1$ or
$\Lambda^3_{27}$. \\

Another distinguished class is $\Lambda^1_7 \oplus \Lambda^3_{27}$. The torsion is given by
\bdm 
T \ = \ - \, * d \omega^3 + *(\beta \wedge \omega^3) \, . 
\edm
Examples of this type were presented and 
discussed in connection with type IIB supergravity solutions in~\cite{GKMW} 
where this expression with a non trivial $\beta$ was first given, though by 
a different approach, making use of both Killing spinor equations.\footnote{We thank J. Gauntlett for drawing our attention to
their work, in particular to the possible class $\Lambda_7^1$ which was missing in a preliminary version of our paper.}. \\

The last class of $G_2$-structures we want to emphasize is $\Lambda^0_1 \oplus \Lambda^1_{7}$. Then $\beta$ is a closed $1$-form and we obtain
\bdm
\delta^g(\omega^3) \ = \ - \, \beta \haken \omega^3, \quad
d \omega^3 \ = \ - \, \lambda \cdot (* \omega^3) + \frac{3}{4} \cdot (\beta 
\wedge \omega^3) \, ,
\edm
\bdm
T \ = \ - \, \frac{\lambda}{6} \cdot \omega^3 + \frac{1}{4} *(\beta \wedge 
\omega^3) \ = \ - \, \frac{7 \lambda}{6} \cdot \omega^3 - \, * d \omega^3 + 
*(\beta \wedge \omega^3) \, .
\edm
In particular, if $ \lambda = 0$, then the torsion form is coclosed,
$\delta^g(T) = 0$.
%

%----------------------------------------------------------------------------
\section{The $\nabla$-Ricci tensor of a $G_2$-structure}
%----------------------------------------------------------------------------
%
We consider a $G_2$-manifold $(M^7,g,\omega^3, \nabla)$ of type 
$\Lambda^0_1 \oplus \Lambda^1_7 \oplus \Lambda^3_{27}$ and its unique
connection $\nabla$ preserving the $G_2$-structure. The totally skew-symmetric
torsion tensor $T$ is a $3$-form of type $\Lambda^3_1 \oplus \Lambda^3_7 
\oplus \Lambda^3_{27}$. For any vector 
field $X$ the covariant derivative $\nabla_XT$ is again a $3$-form of type 
$\Lambda^3_1 \oplus \Lambda^3_7 \oplus \Lambda^3_{27}$. Moreover, 
the $\nabla$-parallel $3$-form $\omega^3$ defines a 
$\nabla$-parallel spinor field $\Psi_0$ (see~\cite{FKMS}). The Clifford 
products $(\nabla_XT) \cdot \Psi_0$ and $(X \haken dT)\cdot \Psi_0$ depend 
only on the $(\Lambda^3_1 \oplus \Lambda^3_7)$-part of the corresponding $3$-forms:
\bdm
\pi^3_1(\nabla_XT) \ := \ \frac{1}{7}\cdot (\nabla_XT, \omega^3) \cdot 
\omega^3, \quad 
\pi^3_7(\nabla_X T) \ := \ \frac{1}{4} \sum \limits^7_{i=1}
(\nabla_XT, e_i \haken *\omega^3) \cdot (e_i \haken * \omega^3) \, ,
\edm
\bdm
\pi^3_1(X \haken dT) \ := \ \frac{1}{7}\cdot(X \haken dT,\omega^3)\cdot \omega^3, 
\quad \pi^3_7(X \haken dT) \ := \ \frac{1}{4} \sum \limits^7_{i=1}
(X \haken dT, e_i \haken *\omega^3) \cdot (e_i \haken * \omega^3) \, .
\edm
The second equation of Corollary 3.2 becomes 
\bdm
\frac{1}{2}(X \haken dT)\cdot\Psi_0 + (\nabla_XT) 
\cdot \Psi_0 -\mbox{Ric}^{\nabla}(X)\cdot\Psi_0 \ = \ 0 \, .
\edm
Using the algebraic formulas
\bdm
(X \haken * \omega^3) \cdot \Psi_0 \ = \  4 \cdot X \cdot \Psi_0, \quad 
\omega^3 \cdot \Psi_0 \ = \ - \, 7 \cdot \Psi_0
\edm
valid for the special spinor $\Psi_0$ related to $\omega^3$ we conclude
\bdm
(X \haken dT, \omega^3) \ = \ - \,2\cdot (\nabla_XT, \omega^3), \quad
\mbox{Ric}^{\nabla}(X) \ = \ \frac{1}{2} \sum \limits^7_{i=1}
(X \haken dT + 2 \cdot \nabla_XT, e_i \haken *\omega^3) \cdot e_i \ .
\edm
The relation between the Ricci tensors 
\bdm
\mbox{Ric}^{\nabla}(X) \ = \ \mbox{Ric}^g(X) + \frac{1}{4} \sum 
\limits^n_{i,j=1} g(T(e_i,X),T(e_j,e_i))\cdot e_j - \frac{1}{2} 
\sum\limits^n_{i=1} \delta^g(T)(X,e_i) \cdot e_i
\edm
allows us to compute the Riemannian Ricci tensor. We summarize all the
derived results in one theorem.\\

{\bf Theorem 5.1.} {\it Let $(M^7,g,\omega^3, \nabla)$ be a $G_2$-manifold 
of type $\Lambda^0_1 \oplus \Lambda^1_7 \oplus \Lambda^3_{27}$ and its unique 
connection $\nabla$ preserving the $G_2$-structure. The Ricci tensor}
$\mbox{Ric}^{\nabla}$ {\it is given by the formula}
\bdm
\mbox{Ric}^{\nabla}(X) \ = \ \frac{1}{2} \sum \limits^7_{i=1}
(X \haken dT + 2 \cdot \nabla_XT, e_i \haken *\omega^3) \cdot e_i \ .
\edm
{\it $T$ is a solution of the equation}
\bdm
\mbox{Ric}^g_{ij} - \frac{1}{4}T_{imn}T_{jmn} - \frac{1}{2}(e_i \haken dT + 
2 \cdot \nabla_{e_i}T, e_j \haken *\omega^3) - \frac{1}{2} \delta^g(T)_{ij} \ = \ 0 \ ,
\edm
{\it and satisfies, for any vector $X$, the condition}
$(X \haken dT, \omega^3) = - \,2 \cdot (\nabla_XT, \omega^3)$. {\it There 
exists a $\nabla$-parallel spinor field $\Psi_0$ such that the $3$-form $\omega^3$ and the Riemannian Dirac operator $D^g$ act on it by} 
\bdm
\omega^3 \cdot \Psi_0 \ = \ - \, 7 \cdot \Psi_0, \quad D^g(\Psi_0) \ = \ 
- \, \frac{3}{4} \cdot T \cdot \Psi_0 \ = \ - \, 
\frac{7}{8} \cdot \lambda \cdot \Psi_0 + \frac{3}{16}\cdot (\beta \haken *\omega^3) \cdot \Psi_0 \, .
\edm
{\bf Example 5.2.} Consider the case of a nearly-parallel $G_2$-structure 
($ \Gamma = \lambda \cdot \mbox{Id}_{T(M^7)}$). 
The torsion form is proportional to the form of the $G_2$-structure, 
\bdm
T \ = \ - \, \frac{\lambda}{6} \cdot \omega^3, \quad d \omega^3 \ = \ 
- \, \lambda \cdot (* \omega^3), \quad dT \ = \ \frac{\lambda^2}{6} \cdot (*\omega^3) \, .
\edm
Then we have
\bdm
\mbox{Ric}^g_{ij}\ = \ \frac{27}{72}\cdot \lambda^2 \cdot \delta_{ij}, \quad  
\frac{1}{4}T_{imn}T_{jmn}\ = \ \frac{3}{72} \cdot \lambda^2 \cdot \delta_{ij}, \quad 
\frac{1}{2}(e_i \haken dT, e_j \haken *\omega^3)\ = \ \frac{24}{72} \cdot 
\lambda^2 \cdot \delta_{ij} \ .
\edm
The formula of Theorem 5.1 for the Ricci tensor generalizes the well-known 
fact that a nearly-parallel $G_2$-manifold is an Einstein space (see~\cite {FKMS}). \\

We consider now only the cocalibrated case, $\beta = 0$. Then the covariant
derivative $\nabla_XT$ is a $3$-form of type $\Lambda^3_1 \oplus \Lambda^3_{27}$ and the formula for the Ricci tensor does not contain the $\nabla_XT$-term.
If $M^7$ is a compact, cocalibrated $G_2$-manifold we can apply the estimate for the first eigenvalue of the Riemannian Dirac 
operator (see~\cite{F0}) 
\bdm
\frac{7}{4 \cdot 6} \| \Psi_0 \|^2 \cdot \mbox{vol}(M^7) \cdot \mbox{Scal}^g_{min} \leq   \int_{M^7}(D^g\Psi_0,D^g\Psi_0) \, ,
\edm
where $\mbox{Scal}^g_{min}$ denotes the minimum of the Riemannian scalar
curvature. Using the equation $8 \cdot D^g(\Psi_0) = - 7 \cdot \lambda \cdot
\Psi_0$ and the definition of the function $\lambda$ we obtain an $L^2$-lower
bound for $(d \omega^3, * \omega^3)$.\\

{\bf Theorem 5.3.} {\it For any compact, cocalibrated $G_2$-manifold the 
following inequality holds}
\bdm
\frac{56}{3} \cdot \mbox{vol}(M^7) \cdot \mbox{Scal}^g_{min} \leq   
\int_{M^7} (d \omega^3, * \omega^3)^2 \ .
\edm
{\it Equality occurs if and only if the cocalibrated $G_2$-structure is nearly
parallel. If the $G_2$-structure is of pure type $\Lambda^3_{27}$, then the minimum of the scalar curvature is non-positive.}\\

The formula of Theorem 5.1 computes in particular the Ricci tensor and the scalar curvature of the connection $\nabla$ :
\bdm
\mbox{Ric}^{\nabla}_{ij} \ = \ \frac{1}{2}(e_i \haken dT, e_j \haken *\omega^3), \quad 
\mbox{Scal}^{\nabla} \ = \ \frac{1}{2} \sum \limits^7_{i=1}
(e_i \haken dT, e_i \haken *\omega^3) \, .
\edm
A cocalibrated $G_2$-structure $(M^7, \omega^3, \nabla, \Psi_0)$ together with its
canonical connection and spinor field solves all of the 
three string equations $\mbox{Ric}^{\nabla}=0, \, \nabla \Psi_0 = 0, \, 
T \cdot \Psi_0 = 0 $ if and only if the $G_2$-structure is geometrically flat ($\lambda = 0 = T$). Therefore we study the first of these equations only.\\ 

{\bf Theorem 5.4.} {\it Let $(M^7,g,\omega^3)$ be a $7$-dimensional Riemannian manifold with a coca\-librated $G_2$-structure. The following conditions are 
equivalent:}
\begin{itemize}
\item[1)]{\it The Ricci tensor} $\mbox{Ric}^{\nabla}$ {\it vanishes.}
\item[2)] {\it The torsion form $T$ is closed and coclosed, $dT = 0 = 
\delta^g(T)$.}
\item[3)] {\it $\lambda$ is constant and the $G_2$-structure} $\omega^3$ {\it satisfies the equation}
\bdm
d*d\omega^3 + \frac{7}{6} \cdot \lambda \cdot d\omega^3 \ = \ 0 \quad 
\lambda \ = \ - \, \frac{1}{7} \cdot (d \omega^3, * \omega^3) \, .
\edm
\end{itemize}
{\it Moreover, in this case we have}
\bdm
\Big(*d\omega^3 + \frac{7}{6} \cdot \lambda \cdot \omega^3 \Big) \wedge d \omega^3 \ = \ 0 \ .
\edm
{\it If the $G_2$-structure is of pure type $\Lambda^0_1$ or
$\Lambda^3_{27}$ and } $\mbox{Ric}^{\nabla}$ {\it vanishes, then the Riemannian manifold $M^7$ is a Ricci flat space with holonomy $G_2$.}\\

{\bf Proof.} The condition $\mbox{Ric}^{\nabla} = 0$ means that $X \haken dT$ is orthogonal to the subspace $\{Y \haken* \omega^3\}= \Lambda^1_7$. 
Moreover, $X \haken dT$ is orthogonal to $\omega^3$ and therefore it belongs to the subspace $\Lambda^3_{27}$. Consequently we obtain for any vector $X$ the conditions
\bdm
(X \haken dT ) \wedge \omega^3 \ = \ 0 \quad \mbox{and} \quad (X \haken dT ) \wedge *\omega^3 \ = \ 0 \, .
\edm
The subspace of all $4$-forms satisfying these algebraic equations is a $G_2$-invariant subspace of $\Lambda^4$ and it is not hard to see that this space is trivial, i.e. we conclude that $dT=0$. Since $\delta^g(T)$ is the antisymmetric part
of the Ricci tensor $\mbox{Ric}^{\nabla}$, the codifferential of the torsion form must vanish, too. These arguments prove the equivalence of the three conditions in the Theorem. The torsion form is 
of type $\Lambda^0_1 \oplus \Lambda^3_{27}$ and we differentiate the equation
$T \wedge \omega^3 = 0$. Then we obtain the last equation of the Theorem.
$\Box$\\

{\bf Remark 5.5.} The cocalibrated $G_2$-structure $(M^7,g,\omega^3,\nabla)$ defines a (homogeneous) solution to the string equations with constant dilation
\bdm
\mbox{Ric}^g_{ij} - \frac{1}{4}T_{imn}T_{jmn} \ = \ 0, \quad \delta^g(T) \ = \ 0  
\edm 
if and only if $\omega^3$ is a solution of the cubic equation
\bdm
d*d\omega^3 + \frac{7}{6} \cdot \lambda \cdot d\omega^3 \ = \ 0 \, ,
\edm
and $7 \cdot \lambda = - \, (d\omega^3,*\omega^3)$ is constant. In this case the torsion form is closed and coclosed. In particular, if $M^7$ is a non-flat $G_2$-structure, we obtain a necessary topological condition:
\bdm
H^3(M^7 ; {\Bbb Z}) \ \neq \ 0 \, .
\edm
We have not succeeded in constructing any cocalibrated $G_2$-structure satisfying this non-linear equation.\\

{\bf Theorem 5.6.} {\it Let $(M^7,g, \omega^3, \nabla)$ be a $7$-dimensional
compact nearly parallel $G_2$-manifold and $\nabla$ be the unique 
$G_2$-connection with totally skew-symmetric torsion. Then every 
$\nabla$-harmonic spinor $\Psi$ is $\nabla$-parallel. Moreover, the space of
$\nabla$-parallel spinors is one-dimensional.}\\

{\bf Proof.} The Dirac operator $D$ is selfadjoint. Let $M^7$ be compact and
consider a $\nabla$-harmonic spinor, $D \Psi =0$. Then Theorem 3.3 implies
\bdm
\sum\limits^7_{i=1} \int_{M^7} (e_i\haken T \cdot\nabla_{e_i}\Psi,\Psi)\ 
= \ 0 \, ,
\edm
since $dT=2 \cdot \sigma^T$ holds in case of a nearly-parallel structure. 
Using the latter equality as well as the Schr\"odinger-Lichnerowicz formula
we obtain
\bdm
 \int_{M^7} \Big(\| \nabla \Psi \|^2 + \frac{1}{2} (dT \cdot \Psi, \Psi)+ \frac{1}{4} 
\mbox{Scal}^{\nabla} \| \Psi \|^2\Big) \ = \ 0 \, .
\edm
The $4$-form $*\omega^3$ acts on spinors as a symmetric endomorphism with the 
eigenvalues $+1$ and $-7$. The result follows now from the estimate
\bdm
\frac{1}{2}(dT \cdot \Psi, \Psi) + \frac{1}{4} \mbox{Scal}^{\nabla} \| 
\Psi \|^2
\geq \Big(\frac{1}{2}\cdot 24 \cdot \lambda^2(-7) + \frac{1}{4} \cdot 48 
\cdot \lambda^2 \cdot 7 \Big) \| \Psi \|^2 \ = \ 0 \, . 
\quad \quad \quad \Box
\edm
{\bf Remark 5.7.} Let us discuss the result from the point of view of the
spectrum of the Riemannian Dirac operator. The first eigenvalue of the 
Riemannian Dirac operator on a compact, simply connected nearly parallel 
manifold $(M^7,g, \omega^3, \nabla)$ is
\bdm
\mu \ = \ \frac{1}{2} \sqrt{\frac{7 \cdot R}{6}} \ = \ \frac{7 \cdot \lambda}
{8} 
\edm
(see \cite{F0}). Here we have constructed a $\nabla$-parallel spinor
\bdm
\nabla \Psi_0 \ = \ \nabla^g_X \Psi_0 - \frac{1}{24} \cdot \lambda \cdot 
(X \haken \omega^3 ) \cdot \Psi_0 \ = \ 0 \, .
\edm
Let us compute the Riemannian Dirac operator :
\bdm
D^g \Psi_0 - \frac{3}{24} \cdot \lambda \cdot \omega^3 \cdot \Psi_0 \ = \ 0 \, .
\edm
The endomorphism $\omega^3$ acts $\Psi_0$ by multiplication by $-7$. 
Therefore we obtain
\bdm
D^g \Psi_0 \ = \ -  \, \frac{7}{8} \cdot \lambda \cdot \Psi_0 \, , 
\edm
i.e. the $\nabla$-parallel spinor field on $M^7$ is the real Killing spinor on $M^7$. In this sense the $\nabla$-parallel spinors on (non-nearly parallel) 
cocalibrated $G_2$-structures generalize the Killing spinors.\\
%
%----------------------------------------------------------------------------
\section{Examples}
%----------------------------------------------------------------------------
%
Denote by $\mbox{H}(3)$ the $6$-dimensional, simply-connected Heisenberg 
group and consider the product $M^7 := \mbox{H}(3) \times {\Bbb R}$. 
There exists a left invariant metric and an orthonormal frame 
$e_1,\, \dots , e_7$ such that the corresponding $3$-form $\omega^3$ 
defines a cocalibrated $G_2$-structure of pure type $\Lambda^3_{27}$ 
(see~\cite{FU}). Indeed, the exterior differentials are given by the 
formulas
\bdm
de_1 \ = \ de_2 \ = \ de_3 \ = \ de_6 \ = \ de_7 \ = \ 0 \ ,
\edm
\bdm
de_4 \ = \ e_1 \wedge e_6 + e_3 \wedge e_7, \quad de_5 \ = \ e_1 \wedge e_3 -
e_6 \wedge e_7 \ 
\edm
and an easy computation yields the following formula for the differential
\bdm
d \omega^3 \ = \ e_1 \wedge e_2 \wedge e_3 \wedge e_4 +  e_2 \wedge e_4 \wedge
 e_6 \wedge e_7 +  e_1 \wedge e_2 \wedge
 e_5 \wedge e_6 -  e_2 \wedge e_3 \wedge
 e_5 \wedge e_7 \ .
\edm 
We see that $d \omega^3 \wedge \omega^3 = 0, \ d \omega^3 \wedge *\omega^3 
= 0$, i.e., the $G_2$-structure is of pure type $\Lambda^3_{27}$. The torsion 
form $T = - \, *d \omega^3$ equals
\bdm
T \ = \ - \, ( e_5 \wedge e_6 \wedge e_7 - e_1 \wedge e_3 \wedge e_5 +
 e_3 \wedge e_4 \wedge e_7 +  e_1 \wedge e_4 \wedge e_6) 
\edm 
and its differential is given by
\bdm
dT \ = \ - \, 4 \cdot e_1 \wedge e_3 \wedge e_6 \wedge e_7 \, .
\edm
The Ricci tensor $2 \cdot \mbox{Ric}^{\nabla}(X,Y) = (X \haken dT, Y \haken 
* \omega^3)$ is a diagonal matrix
\bdm
\mbox{Ric}^{\nabla} \ = \ \mbox{diag}\, (-2,\, 0,\, -2,\, 0,\, 0,\, -2,\, -2)
\edm 
and the scalar curvature becomes negative, $\mbox{Scal}^{\nabla} = - 8$. The
symmetric tensor $T_{imn}T_{jmn}$ is of diagonal form too,
\bdm
T_{imn}T_{jmn} \ = \ \mbox{diag}\, (4,\, 0,\, 4,\, 4,\, 4,\, 4,\, 4) \, ,
\edm
and thus we obtain the Riemannian Ricci tensor
\bdm
\mbox{Ric}^g_{ij} \ = \ \frac{1}{4}T_{imn}T_{jmn} + \mbox{Ric}^{\nabla}_{ij} \ = \ 
\mbox{diag}\, (-1,\, 0,\, -1,\, 1,\, 1,\, -1,\, -1) \, .
\edm
Now we study the $\nabla$-parallel spinors. First of all we need 
the $4$-forms
\bdm
\frac{1}{4}dT + \frac{1}{2} \sigma^T \ = \ - \, e_1 \wedge e_3 
\wedge e_6 \wedge e_7 + (e_3 \wedge e_4 \wedge e_5 \wedge e_6 - 
e_1 \wedge e_4 \wedge e_5 \wedge e_7 - e_1\wedge e_3\wedge e_6\wedge e_7)\, ,
\edm
\bdm
\frac{3}{4}dT - \frac{1}{2} \sigma^T \ = \ -\, 3 e_1 \wedge e_3 
\wedge e_6 \wedge e_7 - \, (e_3 \wedge e_4 \wedge e_5 \wedge e_6 - 
e_1 \wedge e_4 \wedge e_5\wedge e_7 - e_1\wedge e_3\wedge e_6\wedge e_7)\,.
\edm
\newpage
{\bf Lemma 6.1.} 
\begin{itemize}
\item[1)] {\it The $4$-form $dT/4 + \sigma^T/2$ acts in the spinor bundle as a symmetric endomorphism with eigenvalues $(2,\, -4,\, 2,\, 0,\, 2,\,
0,\, 2,\, -4)$.} 
\item[2)] {\it The $4$-form $3 \cdot dT/4 - \sigma^T/2$ acts in the spinor bundle as a symmetric endomorphism with eigenvalues $(2,\, 0,\, 2,\, -4,\, 2,\, - 4,\, 2,\, 0)$.} 
\end{itemize} 

The proof of Lemma 6.1 is an easy computation in the spin representation of the
$7$-dimensional Clifford algebra.\\

{\bf Corollary 6.2.} {\it There are four $\nabla$-parallel spinor fields on 
$M^7$. The torsion form acts trivially on any of these spinors,} $T \cdot \Psi 
= 0$.\\

Let $G$ be a discrete group of isometries acting on $M^7$ and preserving the
$G_2$-structure $\omega^3$. Then $M^7/G$ admits a $G_2$-structure of type
$\Lambda^3_{27}$. \\

{\bf Corollary 6.3.} {\it If $M^7/G$ is a compact manifold and $\Psi$ is
$\nabla$-harmonic, then}
\bdm
6 \cdot \int_{M^7/G} \| \Psi \|^2 \ \geq \ \int_{M^7/G} \| \nabla \Psi \|^2 \ .
\edm 

{\bf Proof.} Combining the Schr\"odinger-Lichnerowicz formula and Theorem 3.3 we obtain, in case of a $\nabla$-harmonic spinor, the equation
\bdm
\int_{M^7/G} \Big(\| \nabla \Psi \|^2 + \frac{1}{4}(dT \cdot \Psi, \Psi) + 
\frac{1}{2}(\sigma^T \cdot \Psi, \Psi) + \frac{1}{4} \mbox{Scal}^{\nabla}
\cdot \|\Psi \|^2 \Big) \ = \ 0 \ .
\edm
Since $\mbox{Scal}^{\nabla} = -8$, the proof follows directly by Lemma 6.1.
 \hfill $\Box$\\

We now discuss a second example. The product $M^7 = N^6 \times {\Bbb R^1}$ of ${\Bbb R^1}$ by a $3$-dimensional complex, solvable Lie group $N^6$ 
admits a left invariant metric such that the following structure equations hold (see~\cite{CF}):
\bdm
de_1 \ = \ 0, \quad de_2 \ = \ 0, \quad de_7 \ = \ 0 \ ,
\edm
\bdm
de_3 \ = \ e_1 \wedge e_3 - e_2 \wedge e_4, \quad de_4 \ = \ e_2 \wedge e_3 
+ e_1 \wedge e_4 \ , 
\edm
\bdm
de_5 \ = \ -\, e_1 \wedge e_5 + e_2 \wedge e_6, \quad de_6 \ = \ -\,  e_2 \wedge e_5 - e_1 \wedge e_6 \ .
\edm 
 
A computation of the exterior products yields the formulas :
\bdm
d*\omega^3 \ = \ 0, \quad d\omega^3 \ = \ 2 \cdot e_1 \wedge e_3 \wedge e_4 
\wedge e_7 - 2 \cdot e_1 \wedge e_5 \wedge e_6 \wedge e_7 \ .
\edm
In particular, the corresponding $G_2$-structure is cocalibrated and of pure type $\Lambda^3_{27}$,
\bdm
*d\omega^3 \wedge \omega^3 \ = \ 0, \quad *d\omega^3 \wedge *\omega^3 \ = \ 0 \ .
\edm

The torsion tensor $T$ of the connection associated with the $G_2$-structure is given by
\bdm
T \ = \ 2 \cdot e_2 \wedge e_5 \wedge e_6 - 2 \cdot 
e_2 \wedge e_3 \wedge e_4, \quad dT \ = \ -4 \cdot e_1 \wedge e_2 \wedge e_5 
\wedge e_6 - 4 \cdot e_1 \wedge e_2 \wedge e_3 \wedge e_4 \ .
\edm

We compute the scalar curvature $\mbox{Scal}^{\nabla} = -16$ of the connection $\nabla$ and the $4$-forms
\bdm
\frac{3}{4}dT - \frac{1}{2} \sigma^T \ = \ -\, 3 \cdot e_1 \wedge e_2 \wedge e_5 \wedge e_6 - 3 \cdot e_1 \wedge e_2 \wedge e_3 \wedge e_4 + 2 \cdot 
e_3 \wedge e_4 \wedge e_5 \wedge e_6 \ ,
\edm
\bdm
\hspace{-0.9cm} \frac{1}{4}dT + \frac{1}{2} \sigma^T \ = \ -\, e_1 \wedge e_2 \wedge e_5 \wedge e_6 - e_1 \wedge e_2 \wedge e_3 \wedge e_4 - 2 \cdot e_3 \wedge e_4 \wedge e_5 \wedge e_6 \ .
\edm
\newpage
{\bf Lemma 6.4.} 
\begin{itemize}
\item[1)] {\it In the spinor bundle the $4$-forms $dT/4 + \sigma^T/2$ acts in the spinor bundle as a symmetric endomorphism with eigenvalues $(4,\, 4,\, -2,\, -2,\, -2,
\, -2,\, 0,\, 0)$.}
\item[2)] {\it The $4$-form $3 \cdot dT/4 - \sigma^T/2$ acts as a symmetric endomorphism with eigenvalues $(4,\, 4,\, 2,\, 2,\, 2,\, 2,\, -8,\, -8)$ in the spinor bundle.}
\end{itemize}

{\bf Corollary 6.5.} {\it There are two $\nabla$-parallel spinor fields on 
$M^7$. The torsion form acts trivially on any of these spinors,} $T \cdot \Psi 
= 0$.\\

Let G be a discrete group of isometries acting on $M^7$ and preserving the
$G_2$-structure $\omega^3$. Then $M^7/G$ admits a $G_2$-structure of type
$\Lambda^3_{27}$. \\

{\bf Corollary 6.6.} {\it If $M^7/G$ is a compact manifold and $\Psi$ is
$\nabla$-harmonic, then}
\bdm
6 \cdot \int_{M^7/G} \| \Psi \|^2 \ \geq \ \int_{M^7/G} \| \nabla \Psi \|^2 \ .
\edm 
\\

We would like to mention that any hypersurface $M^7 \subset {\Bbb R}^8$
admits a cocalibrated $G_2$-structure (see~\cite{FG}). This structure is 
of pure type $\Lambda^0_1$ if and only if the hypersurface is umbilic. The 
pure type $\Lambda^3_{27}$ occurs if and only if the hypersurface is minimal.
The function $\lambda $ is proportional to the mean 
curvature of the hypersurface (see~\cite{FG}). Moreover, it turns out that in 
the decomposition
\bdm
T \ = \ - \, \frac{1}{6} \cdot \lambda \cdot \omega^3 -  \Gamma^3_{27} 
\edm
of the torsion tensor $T$ the $3$-form $\Gamma^3_{27}$ corresponds to the traceless part of the second fundamental form of the hypersurface via the $G_2$-isomorphism $S^2_0({\Bbb R}^7) = \Lambda^3_{27}$. The torsion form is coclosed for hypersurfaces of constant mean curvature. In this case we obtain solutions of the equations in Theorem 5.1 such that $\delta^g(T) = 0$.
%

%----------------------------------------------------------------------------
\section{Sasakian manifolds in dimension five}
%----------------------------------------------------------------------------
The case of the group $G_2$ and dimension $n=7$ discussed in detail fits into a more general approach described in the introduction. We study two further
natural geometric structures: almost metric contact structures and almost hermitian structures. To begin with, let us consider the case of $5$-dimensional 
Sasakian manifolds.\\

{\bf Proposition 7.1.} {\it Every Sasakian manifold $(M^{2k+1},g 
,\xi, \eta, \varphi)$ admits a unique metric connection with totally 
skew-symmetric torsion preserving the Sasakian structure :}
\bdm
\nabla \xi \ = \ \nabla \eta \ = \ \nabla \varphi \ = \ 0 \ . 
\edm
{\it The connection $\nabla$ is given by}
\bdm
g(\nabla_X Y,Z) \ = \ g(\nabla^g_X Y,Z) + \frac{1}{2} (\eta \wedge d \eta)
(X,Y,Z) , \quad T \ = \ \eta \wedge d \eta \ . 
\edm
{\it The torsion form $T$ is $\nabla$-parallel and henceforth
 coclosed,} $\delta^g(T) = 0$. {\it The $4$-form $2 \cdot\sigma^T$ coincides with} $dT$, 
\bdm
2 \cdot \sigma^T \  = \  dT \ = \ d \eta \wedge d \eta \ .
\edm

{\bf Proof.} The existence of the connection has been noticed in earlier papers
(see e.g.\cite{Ko}). The uniqueness will be proved in a more general context in Theorem 8.2. 
\hfill $\Box$\\

We now consider a $5$-dimensional Sasakian manifold $M^5$ and orient it by the condition that the differential of the contact form is given by 
\bdm
d \eta \ = \ 2 \cdot (e_1 \wedge e_2 + e_3 \wedge e_4) \ .
\edm
Furthermore, we fix a spin structure. The endomorphism $\eta \wedge d \eta = 2 \cdot(e_1 \wedge e_2 + e_3 \wedge e_4) \wedge e_5$ acts in the $5$-dimensional spin representation with eigenvalues $(-4, \, 0, \, 0, \, 4)$. Consequently, the spinor bundle splits into two $1$-dimensional and one $2$-dimensional $\nabla$-parallel subbundles. The Clifford multiplication by $\xi$ preserves this decomposition of the spinor bundle and acts on the $1$-dimensional bundles by multiplication by $i$, on the $2$-dimensional bundle by multiplication by $-i$. Suppose that there exists a $\nabla$-parallel spinor $\Psi$. Then one of the subbundles under consideration admits a $\nabla$-parallel spinor 
\bdm
\nabla^g_X \Psi + \frac{1}{4} (X \haken \eta \wedge d\eta ) \cdot \Psi \ = \ 0 
\edm
and the Riemannian Dirac operator for this spinor is given by the formula
\bdm
D^g \Psi + \frac{3}{4} (\eta \wedge d \eta) \cdot \Psi \ = \ 0 \ .
\edm

Let us first discuss the case that $\Psi$ belongs to one of the
$1$-dimensional subbundles defined by the algebraic equation $\xi \cdot  \Psi
= i \cdot \Psi$. In this case we apply Corollary 3.2 
\begin{eqnarray*}
\frac{1}{2} \cdot (X \haken dT)\cdot \Psi - \mbox{Ric}^{\nabla} (X) \cdot
\Psi \ = \ 0 
\end{eqnarray*}

as well as the following algebraic lemma.\\

{\bf Lemma 7.2.} {\it The spinor $\Psi = (1, \, 0, \, 0, \, 0)$ or $(0, \, 0, \, 0, \ 1)$ belongs to the kernel of the endomorphism } 
\bdm
\sum_{1 \leq i<j<k \leq 5} t_{ijk} \cdot e_i \cdot e_j \cdot e_j + \sum_{i=1}^5 x_i \cdot e_i
\edm 
{\it in the $5$-dimensional spin representation if and only if the following equations hold :}
\bdm
x_1 \ = \ - \, t_{234}, \quad x_2 \ = \  t_{134}, \quad x_3 \ = \ - \, t_{124}, \quad 
x_4 \ = \ t_{123}, \quad x_5 \ = \ 0 \ ,
\edm
\bdm
t_{125} \ = \ - \, t_{345}, \quad t_{235} \ = \ - \, t_{145}, \quad 
t_{245} \ = \ t_{135} \ .
\edm

Using these formulas we conclude $\mbox{Ric}^{\nabla}=\mbox{diag}(a, a, a, a, 0)$, where $2 \cdot a:= dT(e_1,e_2,e_3,e_4) = d \eta \wedge d \eta (e_1,e_2,e_3,e_4) = 8$. In particular, we obtain 
\bdm
\mbox{Ric}^g - \mbox{diag}(2,\, 2, \, 2, \, 2, \, 4) \ = \  
\mbox{Ric}^g - \frac{1}{4} T_{imn}T_{jmn} \ = \ \mbox{Ric}^{\nabla} \ ,
\edm
i.e., $\mbox{Ric}^g = \mbox{diag}(6,\, 6, \, 6, \, 6, \, 4) $. A proof similar
to the proof of the existence of Killing spinors on Einstein-Sasakian manifolds (see \cite{FK}) shows that this condition is the only integrability condition for $\nabla$-parallel spinors in the $1$-dimensional subbundles. The endomorphisms act in the spinor bundle 
\bdm
\frac{3}{4} dT - \frac{1}{2} \sigma^T + \frac{1}{4}\mbox{Scal}^{\nabla} \ = \ 4 \cdot (e_1 \cdot e_2 \cdot e_3 \cdot e_4 + 1) \ = \ \frac{1}{4} dT + \frac{1}{2} \sigma^T + \frac{1}{4}\mbox{Scal}^{\nabla} 
\edm
with eigenvalues $0$ and $4$. Corollary 3.2 and Theorem 3.3 yield\\

{\bf Theorem 7.3.} {\it Let $(M^5, g, \xi, \eta, \varphi)$ be a simply connected $5$-dimensional Sasakian spin manifold and consider the unique linear connection $\nabla$ with totally skew-symmetric torsion preserving the Sasakian structure. There exists a $\nabla$-parallel spinor in the subbundle defined by the algebraic equation} $\xi \cdot \Psi = i \cdot \Psi$ {\it if and only if the Riemannian Ricci tensor of $M^5$ has the eigenvalues} $(6, \, 6, \, 6, \, 6, \, 4)$.
{\it A $\nabla$-parallel spinor of this algebraic type is an eigenspinor of the
Riemannian Dirac operator,} $D^g \Psi = \pm \, 3 \cdot \Psi$. {\it In case 
$M^5$ is compact, any $\nabla$-harmonic spinor $\Psi$ is $\nabla$-parallel.} \\

{\bf Example 7.4.} Sasakian manifolds with the described form of the Ricci tensor can be constructed -- for example -- as bundles over $4$-dimensional K\"ahler-Einstein manifolds with positive scalar curvature. Indeed, consider a 
simply connected K\"ahler-Einstein manifold $(N^4,J,g^*)$ with scalar curvature $\mbox{Scal}^* = 32$. Then there exists an $S^1$-bundle $M^5 \rightarrow N^4$ as well as a Sasakian structure on $M^5$ such that the Ricci tensor has the eigenvalues $\mbox{Ric}^g \ = \ \mbox{diag}(6,\, 6, \, 6, \, 6, \, 4)$ 
(see \cite{FKim}). More general, the Tanno deformation of an arbitrary 
$5$-dimensional Einstein-Sasakian structure yields for a special value 
of deformation parameter examples of Sasakian manifolds satisfying the
condition of Theorem 7.3 (see Example 9.3). The Einstein-Sasakian manifolds constructed recently in \cite{BG} admit $\nabla$-parallel spinors with respect
to a Tanno deformation of the Sasakian structure.
\\ 

We discuss the case that the $\nabla$-parallel spinor field $\Psi$ belongs to the $2$-dimensional subbundle defined by the algebraic equation $d \eta \cdot \Psi = 0$. A spinor field of this type is a Riemannian harmonic spinor. Let us again compute the Ricci tensor $\mbox{Ric}^{\nabla}$:\\

{\bf Lemma 7.5.} {\it The spinor $\Psi = (0, \, 1, \, 0, \, 0)$ belongs 
to the kernel of the endomorphism } 
\bdm
\sum_{1 \leq i<j<k \leq 5} t_{ijk} \cdot e_i \cdot e_j \cdot e_j + \sum_{i=1}^5 x_i \cdot e_i
\edm 
{\it in the $5$-dimensional spin representation if and only if the following equations hold :}
\bdm
x_1 \ = \ t_{234}, \quad x_2 \ = \ - \, t_{134}, \quad x_3 \ = \ t_{124}, \quad 
x_4 \ = \ - \, t_{123}, \quad x_5 \ = \ 0 \ ,
\edm
\bdm
t_{125} \ = \ t_{345}, \quad t_{235} \ = \ t_{145}, \quad 
t_{245} \ = \ - \, t_{135} \ .
\edm

In this case we obtain 
\bdm
\mbox{Ric}^{\nabla} \ = \ \mbox{diag}(-4, \, -4, \, -4, \, -4, \, 0), 
 \quad \mbox{Ric}^g \ = \ \mbox{diag}(-2,\, -2, \, -2, \, -2, \, 4) \ .  
\edm

We compute the endomorphisms acting in the spinor bundle :
\bdm
\frac{3}{4} dT - \frac{1}{2} \sigma^T + \frac{1}{4}\mbox{Scal}^{\nabla} \ = \ 4 \cdot (e_1 \cdot e_2 \cdot e_3 \cdot e_4 - 1) \ = \ \frac{1}{4} dT + \frac{1}{2} \sigma^T + \frac{1}{4}\mbox{Scal}^{\nabla} \ .
\edm
The Clifford product $e_1 \cdot e_2 \cdot e_3 \cdot e_4$ acts in the $2$-dimensional subbundle of the spin bundle as the identity. Corollary 3.2 and Theorem 3.3 yield\\ 

{\bf Theorem 7.6.} {\it Let} $(M^5, g, \xi, \eta, \varphi)$ {\it be a 5-dimensional Sasakian spin manifold and consider the unique linear connection $\nabla$ with totally skew-symmetric torsion preserving the Sasakian structure. If there exists a $\nabla$-parallel spinor in the subbundle defined by the algebraic equation} $d\eta \cdot \Psi = 0$, {\it then the Riemannian Ricci tensor of $M^5$ has the eigenvalues $(-2, \, -2, \, -2, \, -2, \, 4)$. Any $\nabla$-parallel 
spinor in this $2$-dimensional subbundle satisfies the equations}
\bdm
\nabla^g_{\xi}\Psi \ = \ 0, \quad \nabla^g_{X}\Psi \ = \ \frac{1}{2} \, 
\varphi(X) \cdot \xi \cdot \Psi \ = \ -  \frac{i}{2} \, \varphi(X) \cdot \Psi, \quad 
d\eta \cdot \Psi \ = \ 0 \, .
\edm
{\it In particular, it is harmonic with respect to the Riemannian connection. Any $\nabla$-harmonic spinor $\Psi$ on a compact manifold $M^5$ satisfying the 
algebraic condition $d \eta \cdot \Psi = 0$ is $\nabla$-parallel.} \\

{\bf Example 7.7.} In ${\Bbb R}^5$ we consider the $1$-forms
\bdm
e_1 \ = \ \frac{1}{2}\cdot dx_1, \quad e_2 \ = \ \frac{1}{2}\cdot dy_1, \quad 
e_3 \ = \ \frac{1}{2}\cdot dx_2, \quad e_4 \ = \ \frac{1}{2}\cdot dy_2, \, ,
\edm
\bdm
e_5 \ = \ \eta \ = \frac{1}{2}\cdot(dz - y_1 \cdot dx_1 - y_2 \cdot dx_2) \, .
\edm
We obtain a Sasakian manifold (see \cite{Blair}) and it is not hard to see that it admits $\nabla$-parallel spinors of type $d \eta \cdot \Psi = 0$. The 
Sasakian structure arises from left invariant vector fields on a $5$-dimensional Heisenberg group.\\

Sasakian manifolds with $\nabla$-parallel spinors of
type $F \cdot \Psi = 0$ may be constructed as bundles over the $4$-dimensional 
torus. Indeed, suppose that 
the Sasakian structure is regular. Then $M^5$ is a $S^1$-bundle over $N^4$. 
The spinor field $\Psi$ is projectable and induces a parallel spinor field 
$\Psi^*$ in the negative spinor bundle $\Sigma^-(N^4)$ over $N^4$ (see \cite{Mor}). Consequently, $N^4$ is a selfdual, Ricci-flat K\"ahler manifold. On the 
other hand, the endomorphism $\varphi$ projects too and we obtain a second 
integrable, but in general not parallel positive complex structure. 
There is only one possibility for $N^4$, the torus $T^4$. 
In a forthcoming paper we will study these $\nabla$-parallel spinor even for 
normal almost contact metric structures in more details. In particular it
turns out that the Example 7.7 is (locally) the only Sasakian space with
$\nabla$-parallel spinors of type $d \eta \cdot \Psi = 0$.
%
%----------------------------------------------------------------------------
\section{Almost contact connections with totally skew-symmet\-ric torsion}
%----------------------------------------------------------------------------
Let us discuss the latter results from a more general point of view and 
consider an almost contact metric manifold $(M^{2k+1},g,\xi,\eta,\varphi)$, 
i.e., a Riemannian manifold equipped with a 1-form $\eta$, a (1,1)-tensor $\varphi$ and a vector field $\xi$ dual to $\eta$ with respect to the metric $g$ such that the following compatibility conditions are satisfied (see \cite{Blair}):
\bdm
\eta(\xi) \ = \ 1,\quad \varphi^2 \ = \ - \, \mbox{Id} + \eta\otimes\xi, \quad
g(\varphi (X),\varphi(Y) ) \ = \ g(X,Y) - \eta(X) \cdot \eta(Y), \quad \varphi (\xi) \ = \ 0.
\edm
 
Let us introduce the fundamental form $F(X,Y) :=g(X,\varphi (Y))$ as well as the Nijenhuis tensor
$$ N(X,Y) \ := \ [\varphi (X),\varphi (Y)] + \varphi^2[X,Y] -\varphi[\varphi (X),Y] - \varphi[X,\varphi (Y)] + d\eta(X,Y) \cdot \xi \ ,$$
$$
N^2(X,Y) \ := \ d\eta(\varphi(X),Y) + d\eta(X,\varphi(Y))\ .
$$ 
We recall some notions describing the different types of almost contact metric structures. If $2 \cdot F=d\eta$, then we have a {\it contact metric structure}, if $N=0$, we have a {\it normal contact structure}. A {\it K-contact structure} is a 
 contact metric structure such that the vector field $\xi$ is a Killing vector 
field.  If the structure is normal and K-contact, then it is a {\it Sasaki 
structure} (a complete classification of almost contact metric structures is 
presented in \cite{AG,CG,CM}). The Nijenhuis tensor of type $(0, \, 3)$ is given by $N(X,Y,Z)=g(N(X,Y),
Z)$. We have the following general identities \cite{Blair,CG,CM}:
\begin{eqnarray*}
2 \cdot g((\nabla^g_X\varphi)(Y),Z) &=& dF(X,\varphi(Y),\varphi(Z)) - dF(X,Y,Z) +
N(Y,Z,\varphi(X)) \\
&& + \eta(X) \cdot N^2(Y,Z) + \eta(Z) \cdot d \eta(\varphi(Y),X) +
\eta(Y) \cdot d \eta(X,\varphi(Z)) \ ,
\end{eqnarray*}
\begin{eqnarray*}
g((\nabla^g_X\varphi)(Y),Z) + g((\nabla^g_X\varphi)(\varphi(Y)),\varphi(Z))
= \eta(Y) \cdot (\nabla^g_X\eta)(\varphi(Z))- \eta(Z) \cdot (\nabla^g_X\eta)(\varphi(Y)) \ ,
\end{eqnarray*}
\bdm
g((\nabla^g_X\varphi)(\varphi(Y)),\xi) \ = \ (\nabla^g_X\eta)(Y) \ = \ g(\nabla^g_X \xi,Y) \ ,
\edm
\begin{eqnarray*}
N(X,Y,Z) &=& - N(\varphi (X),\varphi(Y),Z) + \eta(X) \cdot N(\xi,Y,Z) + \eta(Y)\cdot N(X,\xi,Z)\\
&=& - N(\varphi (X),Y,\varphi(Z)) + \eta(Z)\cdot N(\xi,X,Y) - \eta(X) \cdot 
N(\xi,\varphi(Y),\varphi(Z)) \ .
\end{eqnarray*}
Finally, we introduce the forms
\begin{eqnarray*}
dF^-(X,Y,Z) &:=& dF(X,\varphi(Y),\varphi(Z)) + dF(\varphi(X),Y,\varphi(Z)) \\ 
& &  + \, dF(\varphi(X),\varphi(Y),Z)  - dF(X,Y,Z) \, , \\
d^{\varphi}F(X,Y,Z) &:=& - dF(\varphi(X),\varphi(Y),\varphi(Z)) \, ,
\end{eqnarray*}

and a direct consequence of the definitions is the following \\

{\bf Proposition 8.1.} On any almost contact metric manifold the identities 
hold:
\begin{itemize}
\item[1)] $dF^-(X,Y,Z) =  - N(X,Y,\varphi(Z)) - N(Y,Z,\varphi(X)) - N(Z,X,\varphi(Y)) \,$ ,
\item[2)] $N(X,Y) =  (\nabla^g_{\varphi(X)}\varphi)(Y) -  (\nabla^g_{\varphi(Y)}\varphi)(X) +  (\nabla^g_X\varphi)(\varphi(Y)) - (\nabla^g_Y\varphi)
(\varphi(X)) $

$\quad\quad\quad\quad\quad - \eta(Y) \cdot \nabla^g_X\xi + \eta(X) \cdot 
\nabla^g_Y\xi\,.$
\end{itemize}

A linear connection $\nabla$ is said to be an almost contact connection if it
preserves the almost contact structure :
\bdm
\nabla g \ = \ \nabla \eta \ = \ \nabla \varphi \ = \ 0 \ .
\edm

{\bf Theorem 8.2.}
{\it Let $(M^{2k+1},g,\xi,\eta,\varphi)$ be an almost contact metric manifold. 
The following conditions are equivalent:
\begin{itemize}
\item[$1$)] The Nijenhuis tensor $N$ is skew-symmetric and $\xi$ is a Killing vector field.
\item[$2$)] There exists an almost contact linear connection $\nabla$ with 
totally skew-symmetric torsion tensor $T$. 
\end{itemize}
Moreover, this connection is unique and determined by }
\bdm
g(\nabla_XY,Z) \ = \ g(\nabla^g_XY,Z) + \frac{1}{2}T(X,Y,Z) \ ,
\edm
{\it where $\nabla^g$ is the Levi-Civita connection and the torsion $T$ is defined by} 
\bdm
T \ = \ \eta \wedge d\eta + d^{\varphi}F + N - \eta \wedge ( \xi \haken N) \ .
\edm

{\bf Proof.} Let us assume that such a connection exists. Then
\bdm
0 \ = \ g(\nabla^g_X\xi,Z) + \frac{1}{2}T(X,\xi,Z) 
\edm
holds and the skew-symmetry of $T$ yields that $\xi$ is a Killing vector 
field, $d\eta = \xi \haken T, \ \xi \haken d\eta=0$ and 
\bdm
T(\varphi(X),\varphi(Y),Z) - T(X,Y,Z) +T(\varphi(X),Y,\varphi(Z)) + 
T(X,\varphi(Y),\varphi(Z)) \  = \  - \, N(X,Y,Z) \, .
\edm
The latter formula shows that $N$ is skew-symmetric. Since $\varphi$ is
$\nabla$-parallel, we can express the Riemannian covariant derivative of $\varphi$ by the torsion form:
\bdm
T(X,Y,\varphi(Z)) + T(X,\varphi(Y),Z) \ = \ - \, 
2 \cdot g\big((\nabla^g_X\varphi)Y,Z \big) \, .
\edm
Taking the cyclic sum in the above equality, we obtain
\bdm
\sigma_{X,Y,Z} \, T(X,Y,\varphi(Z)) \ = \ - \, \sigma_{X,Y,Z} \, g((\nabla^g_X\varphi)
Y,Z) \, .
\edm
We use Proposition 8.1 as well as the identity preceding it to get
\bdm
- \, \sigma_{X,Y,Z} \, T(X,Y,\varphi(Z)) \ = \ 
\sigma_{X,Y,Z} \, g((\nabla^g_X\varphi)Y,Z) \ = \  - \, dF(X,Y,Z) \, .
\edm
Adding this result to the formula expressing the Nijenhuis tensor $N$ by 
the torsion $T$, some calculations yield 
\bdm
T(\varphi(X),\varphi(Y),\varphi(Z)) \ = \ dF(X,Y,Z) -N(X,Y,\varphi(Z)) - 
\eta(Z) \cdot N^2(X,Y).
\edm
By replacing $X,Y,Z$ by $\varphi(X),\varphi(Y),\varphi(Z)$ and using the 
symmetry property of the Nijenhuis tensor mentioned before 
Proposition 8.1, we obtain the formula for the torsion tensor $T$. \\
For the converse, suppose that the almost contact structure has the 
properties 1) and define the connection $\nabla$ by the formulas in 2). 
Clearly $T$ is 
skew-symmetric and $\xi \haken T=d\eta=2\nabla^g\eta$. Since $\xi$ is a 
Killing vector field, we conclude $\nabla g=\nabla \xi =0$. Furthermore, 
using the conditions 1) and Proposition 8.1, we obtain $\xi \haken dF = N^2$. Finally we have to prove that
$\nabla \varphi =0$. This follows by straightforward computations using the 
relation between $\nabla \varphi$ and the torsion tensor $T$, 
Proposition 8.1 as well as the following lemma. \hfill $\Box$\\

{\bf Lemma 8.3.} {\it Let $(M^{2k+1},g,\xi,\eta,\varphi)$ be an almost contact
metric manifold with a totally skew-symmetric Nijenhuis tensor $N$ and Killing vector field $\xi$. Then the 
following equalities hold:}
\bdm
\nabla^g_{\xi}\xi \ = \ \xi \haken d \eta \ = \ 0 \ ,   
\edm
\bdm
N(\varphi(X),Y,\xi) \ = \ N(X,\varphi(Y),\xi) \ = \ N^2(X,Y) \ = \ 
dF(X,Y,\xi) \ = \ - \, dF(\varphi(X),\varphi(Y),\xi) \ .
\edm

{\bf Proof.} The identities in and before Proposition 8.1 imply
\bdm
0 \ = \ N(X,\xi,\xi) \ = \ (\nabla^g_{\xi}\varphi)(\xi), \quad
0\ = \ N(\xi,X,X) \ = \ (\nabla^g_X\eta)(X)-(\nabla^g_{\varphi(X)}\eta)(\varphi (X)) \ ,
\edm
\bdm
N(X,Y,\xi) \ = \ d\eta(X,Y)-d\eta(\varphi(X),\varphi(Y)) \ .
\edm
Hence, the first two equalities follow. We take the cyclic sum in the second idendity before Proposition 8.1 and put $X=\xi$ to obtain
\bdm
- \, dF(\xi,Y,Z)+g\big((\nabla^g_{\xi}\varphi)(\varphi(Y)),\varphi(Z)\big) \ = \ - \, 
(\nabla^g_Y\eta)(\varphi(Z)) + (\nabla^g_Z\eta)(\varphi(Y)) \ .
\edm
On the other hand, using again the general formulas we calculate that
\bdm
- \, dF(\xi,\varphi(Y),\varphi(Z)) \ = \ g((\nabla^g_{\xi}\varphi)(\varphi(Y)),\varphi(Z))+(\nabla^g_{\varphi(Y)}\eta)(Z) - (\nabla^g_{\varphi(Z)}\eta)(Y) \ .
\edm
Summing up the latter two equalities we obtain the last equalities in the 
lemma since $\xi$ is a Killing vector field. \hfill $\Box$\\

We discuss these results for some special contact structures.\\

{\bf Theorem 8.4.} {\it Let $(M^{2k+1},g,\xi,\eta,\varphi)$ be an almost 
contact metric manifold with totally skew-symmetric Nijenhuis tensor $N$. 
Then the condition $dF=0$ implies $N=0$.} 
\begin{itemize}
\item[1)] {\it A contact metric structure ($2 \cdot F=d\eta$) admits an almost contact
connection with totally skew-symmetric torsion if and only if it is Sasakian. 
In this case, the connection is unique, its torsion is given by}
\bdm
T \ = \ \eta \wedge d \eta
\edm
{\it and $T$ is parallel, $\nabla T = 0$.}
\item[2)] {\it A normal ($N=0$) almost contact structure admits a unique 
almost contact connection with totally skew-symmetric torsion if and only if $\xi$ is a Killing vector field. The torsion $T$ is then given by} 
\bdm
T \ = \ \eta \wedge d\eta  + d^{\varphi}F.
\edm
\end{itemize}
{\bf Proof.} If $dF=0$, Lemma 8.3 implies that $N^2= \xi \haken N=0$. Then
Proposition 8.1 leads to $0=dF^-(X,Y,Z)=-3 \cdot N(\varphi(X),Y,Z)$. The
assertion that $\nabla T=0$ in a Sasakian manifold follows by direct
verification. \hfill $\Box$\\
\section{Almost contact structures, parallel spinors and holonomy  
group}

Let $(M^{2k+1},g,\xi,\eta,\varphi)$ be a $(2k+1)$-dimensional almost
contact metric manifold with totally skew-symmetric Nijenhuis
tensor $N$ and Killing vector $\xi$ and denote by $\nabla$ the unique
almost contact connection with a totally skew-symmetric torsion (Theorem 8.2).
 Since $\nabla\xi=0$ the (restricted) holonomy group $\mbox{Hol}^{\nabla}$ of 
$\nabla$ is contained in $U(k)$. This group cannot occur as the isotropy
group of any spinor. The spinor bundle $\Sigma$ of a contact spin
manifold decomposes under the action of the fundamental form $F$ into the sum 
(see~\cite{FKim})
\bdm
\Sigma \ = \ \Sigma_0 \oplus \ldots \Sigma_k , \quad \mbox{dim}(\Sigma_r) \ 
= \ \left(\!\!\begin{array}{c}k \\ r\end{array}\!\! \right) \ . 
\edm
The isotropy group of a spinor of type $\Sigma_0$ or $\Sigma_k$ coincides with the subgroup $SU(k) \subset U(k)$. Consequently, there exists locally a $\nabla$-parallel spinor of type $\Sigma_0, \Sigma_k$ if and only if  $\mbox{Hol}^{\nabla}$ is contained in $SU(k)$. We shall express this condition in terms of the curvature of $\nabla$. The group $\mbox{Hol}^{\nabla}$ is contained in $SU(k)$ if
and only if the $2$-form
\bdm
\varrho^{\nabla}(X,Y) \ := \ \frac{1}{2}\sum_{i=1}^{2k+1}R^{\nabla}(X,Y,e_i,
\varphi(e_i)) \ .
\edm
vanishes,  $\varrho^{\nabla}=0$. Let us introduce the torsion 1-form
$\omega^{\nabla}$ as well as the 2-form $\lambda^{\nabla}$ by
\bdm
\omega^{\nabla}(X) \ := \ -\frac{1}{2}\sum_{i=1}^{2k+1}T(X,e_i,\varphi
(e_i))\, , \quad 
\lambda^{\nabla}(X,Y) \ := \ \frac{1}{2}\sum_{i=1}^{2k+1}dT(X,Y,e_i,\varphi
(e_i)) \ .
\edm
{\bf Proposition 9.1.} {\it Let $(M^{2k+1},g,\xi,\eta,\varphi)$ be a 
$(2k+1)$-dimensional almost contact metric manifold with totally 
skew-symmetric Nijenhuis tensor $N$ and Killing vector $\xi$. Let $\nabla$ 
be the unique almost contact connection with totally skew-symmetric torsion. 
Then one has} 
\bdm
\varrho^{\nabla}(X,Y) \ = \ \mbox{Ric}^{\nabla}(X,\varphi(Y))-(\nabla_X\omega^
{\nabla})(Y)+\frac{1}{4}\lambda^{\nabla}(X,Y) \ .
\edm

{\bf Proof.} We follow the scheme in \cite{IP}, Section 3, and use the
curvature properties of $\nabla$ from Section 2 to calculate
$\lambda^{\nabla}(X,Y)$ : 
\bdm
-2(\nabla_X\omega)(Y)+2(\nabla_Y\omega)(X) +
2\sum_{i=1}^{2k+1}\Big( \sigma^T(X,Y,e_i,\varphi(e_i))-
(\nabla_{\varphi(e_i)}T)(X,Y,e_i) \Big) \ .
\edm
The first Bianchi identity for $\nabla$ together with the latter
identity implies 
\bdm
4\rho^{\nabla}(X,Y)+2\mbox{Ric}^{\nabla}(Y,\varphi(X))-2\mbox{Ric}^{\nabla}
(X,\varphi(Y)) \ = \ \lambda^{\nabla}(X,Y)-2(\nabla_X\omega)(Y)+
2(\nabla_Y\omega)(X) \ .
\edm
Using the relation between the curvature tensors of $\nabla$ and
$\nabla^g$, we obtain
\bdm
\mbox{Ric}^{\nabla}(Y,\varphi(X))+\mbox{Ric}^{\nabla}(X,\varphi(Y)) \ 
= \ (\nabla_X\omega^{\nabla})(Y)+(\nabla_Y\omega^{\nabla})(X) \ .
\edm
The last two equalities lead to the desired formula.  $\hfill \Box$\\

We apply Proposition 9.1 in case of a Sasakian manifold.\\

{\bf Theorem 9.2.} {\it Let $(M^{2k+1},g,\xi,\eta,\varphi,\nabla)$ be a simply connected $(2k+1)$-dimensional Sasakian spin manifold and $\nabla$ be the 
unique almost contact connection with totally skew-symmetric torsion.
Then there exists a $\nabla$-parallel spinor of type $\Sigma_0$ or $\Sigma_k$
if and only if the Ricci tensor is given by the formula :}
\bdm
\mbox{Ric}^{\nabla}\ = \ 4 \cdot (k-1) \cdot (g-\eta\otimes \eta) \ . 
\edm
{\it This condition is equivalent to}
\bdm
\mbox{Ric}^g \ = \ 2 \cdot (2k-1) \cdot g- 2 \cdot (k-1) \cdot \eta\otimes 
\eta \ .
\edm

{\bf Proof.} On a Sasakian manifold $T=\eta\wedge d\eta= 2 \cdot \eta\wedge
F$ and $\nabla T=0$, where $F(X,Y) = g(X, \varphi(Y))$ is the fundamental form of the Sasakian structure. Consequently, we calculate that
\bdm
\nabla(\omega^{\nabla}) \ = \ 0,\quad \lambda^{\nabla} \ = \ 16 \cdot 
(1-k)F,\quad \sum_{i=1}^{2k+1}g(T(X,e_i),T(Y,e_i)) \ = \ 8 \cdot g+ 8 \cdot 
(k-1)\eta\otimes \eta \ ,
\edm
and the proof follows from Proposition 9.1.  $\hfill \Box$\\

{\bf Remark 9.3.} Sasakian manifolds with the prescribed form of the Ricci 
tensor admit Sasakian quasi-Killing spinors of type $(\pm 1/2,b)$
(see~\cite{FKim}, Theorem 6.3). The Tanno deformation of an Einstein-Sasakian 
manifold defined by the formulas
\bdm
\tilde{\varphi} \ := \ \varphi, \quad \tilde{\xi} \ := \ a^2 \cdot \xi,
\quad \tilde{\eta} \ := \ a^{-2} \cdot \eta, \quad \tilde{g} \ := \ 
a^{-2} \cdot g + (a^{-4}-a^{-2}) \cdot \eta \otimes \eta 
\edm
yields for the parameter $a^2:= 2k/(k+1)$ a Sasakian manifold satisfying
the condition of Theorem 9.2 and vice versa (see~\cite{FKim}, Lemma 6.7 and Lemma 6.8).
%
%

%----------------------------------------------------------------------------
\section{Almost hermitian connections with totally skew-symmet\-ric torsion}
%----------------------------------------------------------------------------
%
In this section we study connections with totally skew-symmetric torsion and 
preserving an almost complex structure. These exist, for example, for nearly K\"ahler manifolds. In dimension $n=6$ nearly K\"ahler manifolds have special properties and they are precisely the 6-dimensional manifolds admitting real Killing spinors (see~\cite{Gru1} and~\cite{Gru2}).\\ 

{\bf Theorem 10.1.} {\it Let $(M^{2n},g, J)$ be a $2n$-dimensional almost 
complex manifold. Then there exists a linear connection with totally 
skew-symmetric torsion preserving the hermitian structure $(g, J)$ if and only if the Nijenhuis tensor $N(X,Y,Z):=g(N(X,Y),Z)$ is a $3$-form. In this case the 
connection is unique and is determined by
\bdm
T(X,Y,Z) \ = - \, d \Omega (J(X), J(Y), J(Z)) + N(X,Y,Z) \, ,  
\edm
where $\Omega$ is the K\"ahler form.}\\

{\bf Proof.} The result can be derived from the considerations in
\cite{Gau1} or from Proposition 4.1. We sketch a direct proof. 
Since $\nabla g= \nabla J=0$, we have
$$ T(J(X), J(Y), Z) - T (X,Y,Z) + T( J(X), Y, J(Z))+ T(X,J(Y), J(Z)) \ = \ - \, N( X,Y,Z) \, ,$$
which shows that $N$ is a $3$-form. The formula for the torsion form follows 
from the following identities on an almost complex manifold with skew 
symmetric tensor $N$.
\begin{eqnarray*}
2 \cdot g \Big( (\nabla_X^g J)Y,Z \Big) &=& d \Omega (X, J(Y), J(Z)) - 
d \Omega (X,Y,Z)+N(Y,Z, J(X))\\
4 \cdot d \Omega^- (X,Y,Z) &:=& d \Omega (X, J(Y), J(Z)) - d \Omega (X,Y,Z)+
d \Omega (J(X), J(Y), Z)\\
&+& d \Omega (J(X), Y, J(Z))\ 
= \ - \, 3 \cdot N(J(X), Y, Z) \ .  \mbox{} \quad \quad \quad \quad \quad \quad \quad  \hfill \Box
\end{eqnarray*} 

{\bf Corollary 10.2.} {\it On an almost K\"ahler manifold there does not exist
 a hermitian connection with totally skew-symmetric torsion.}\\

{\bf Corollary 10.3.} {\it On any nearly K\"ahler manifold the torsion form $T$ is $\nabla$-parallel and henceforth coclosed, $\delta^g(T) = 0$.} \\

{\bf Proof.} If $(M^{2n},g,J)$ is a nearly K\"ahler manifold, then 
\bdm
d\Omega(X,Y,Z) \ = \ - \, d\Omega^-(X,Y,Z) \ = \ 3 \cdot N(J(X),Y,Z)
\edm

and Theorem 10.1 yields that $4 \cdot T= N$. Moreover, $\nabla$ is the 
characteristic connection considered by Gray \cite{Gra} and $T$ is $\nabla$-parallel, $\nabla T= \nabla N=0$ (see e.g. \cite{Kir,MoS}). $\quad  \hfill \Box$ \\

We compute the Ricci tensor $\mbox{Ric}^{\nabla}$ for a $6$-dimensional nearly K\"ahler spin manifold.\\

{\bf Proposition 10.4.} {\it On a $6$-dimensional nearly K\"ahler manifold with 
non-vanishing Nijenhuis tensor $N \neq 0$ we have}
\bdm
T_{imn}T_{jmn} \ = \ 2 \cdot a \cdot g_{ij}, \quad \mbox{Ric}^g \ = \ \frac{5}{2} 
\cdot a \cdot g, \quad \mbox{Ric}^{\nabla} \ = \ 2 \cdot a \cdot g \ .
\edm
{\it The $4$-form $2 \cdot\sigma^T$ coincides with} $dT$, 
\bdm
2 \cdot \sigma^T \ = \ dT \ = \ a \cdot (\Omega \wedge  \Omega) \ .
\edm

{\bf Proof.} We recall (see \cite{Gra}) that any $6$-dimensional nearly 
K\"ahler manifold is Einstein and of constant type, i.e.
\bdm
\mbox{Ric}^g \ = \ \frac{5}{2} \cdot a \cdot g, \quad 
||(\nabla^g_XJ)Y||^2 \ = \ \frac{1}{2} \cdot a \cdot \big(||X||^2 \cdot
||Y||^2 - g^2(X,Y) - g^2(X,J(Y)) \big)
\edm
where $a := \mbox{Scal}^g/15$ is a positive constant. Polarizing the latter
equality und using the identity $ 4 \cdot J(\nabla^g_{X}J)Y = N(X,Y) = 
4 \cdot T(X,Y)$ we get
\bdm
T_{imn}T_{jmn} \ = \ 2 \cdot a \cdot g_{ij}, \quad 2 \cdot \sigma^T \ = \
a \cdot \Omega \wedge \Omega \ .
\edm

We calculate $\mbox{Ric}^{\nabla}_{ij} = \mbox{Ric}^g_{ij} - \frac{1}{4} \cdot
 T_{imn}T_{jmn} = 2 \cdot a \cdot g_{ij}$ and the result follows.  $\quad  \hfill \Box$ \\

We consider again the general almost complex case. Let $(M^{2n},g,J)$ be an
almost complex manifold with totally skew-symmetric
tensor $N$. Then $M^{2n}$ is of type $G_1$ according to Gray-Hervella
classification (see~\cite{GH}). Denote by $\nabla$ the unique hermitian 
connection with totally skew-symmetric torsion $T$ described in Theorem 9.1. 
The Ricci form of $\nabla$ is defined by
\bdm
\varrho^{\nabla}(X,Y) \ = \ \frac{1}{2}\sum_{i=1}^{2n} R^{\nabla}(X,\, Y,\, 
e_i,\, J(e_i)) \ . 
\edm
The holonomy group $\mbox{Hol}^{\nabla}$ of $\nabla$ is contained in 
$SU(n)$ if and only if $\varrho^{\nabla}=0$. We define the Lee form $\theta$ and 
the tensor $\lambda^{\omega}$ by
\bdm
\theta(X) \ = \ - \, \frac{1}{2} \sum_{i=1}^{2n} T(J(X),\, e_i,\, J(e_i)), 
\quad \lambda^{\omega}(X,Y) \ = \ \sum_{i=1}^{2n} dT(X,\, Y,\, e_i,\, J(e_i)) \ .
\edm 
We remark that the formula (3.16) in \cite{IP} holds in the
general case of a $G_1$-manifold, 
\bdm
\varrho^{\nabla}(X,Y) \ = \ \mbox{Ric}^{\nabla}(X,J(Y))+(\nabla_X\theta)J(Y)
+ \frac{1}{4} \lambda^{\omega}(X,Y) \ .
\edm
{\bf Theorem 10.5.} {\it Let $(M^{2n},g,J,\nabla)$ be an almost hermitian
manifold of type $G_1$ with its unique linear connection $\nabla$ with
totally skew-symmetric torsion. Then $\mbox{Hol}^{\nabla} \subset SU(n)$
if and only if }
\bdm
0 \ = \ \mbox{Ric}^{\nabla}(X,J(Y))+(\nabla_X\theta)J(Y)
+ \frac{1}{4} \lambda^{\omega}(X,Y) \ .
\edm
{\bf Corollary 10.6.} {\it Let $(M^{2n},g,J,\nabla)$ be a nearly K\"ahler 
manifold. Then} $\mbox{Hol}^{\nabla}$ {\it is contained in SU(n) if and only 
if} 
\bdm
\mbox{Ric}^{\nabla}(X,Y) \ = \ \frac{1}{4} \lambda^{\omega}(X,J(Y)) \ .
\edm
{\it If the manifold is not K\"ahler and the dimension is $6$, then the above condition is always satisfied.}\\                     

{\bf Proof.} In the nearly K\"ahler case the torsion form $T$ is $\nabla$-parallel and therefore $\theta$ is parallel too, $\nabla
\theta=0$. In the 6-dimensional strictly
nearly K\"ahler case the condition of the Corollary $10.6$ is a
consequence of Proposition $10.4.$  $\quad  \hfill \Box$\\

Finally we study in detail the $\nabla$-parallel and $\nabla$-harmonic spinor
fields on a $6$-dimensional nearly K\"ahler manifold $M^6$. The spinor bundle $\Sigma(M^6)$ splits into the $1$-dimensional subbundles $E^{\pm} \subset \Sigma^{\pm}(M^6)$ of the spinor bundle defined by the equation
\bdm
J(X) \cdot \Psi^{\pm} \ = \ \mp \, i \cdot X \cdot \Psi^{\pm}
\edm
as well as their orthogonal complements. The connection $\nabla$ preserves this decomposition. We discuss the 
integrability condition (see Corollary 3.2) for the existence of a 
$\nabla$-parallel spinor. A purely algebraic computation in the $6$-dimensional spin representation proves the following\\

{\bf Lemma 10.7.} {\it A spinor $\Psi^{\pm}$ satisfies the equation}
\begin{eqnarray*}
\frac{1}{2} (X \haken dT) \cdot \Psi^{\pm} - \mbox{Ric}^{\nabla} (X) \cdot
\Psi^{\pm} \ = \ 0 
\end{eqnarray*}
{\it for all vectors $X$ if and only if it belongs to $E^{\pm}$. The endomorphism}
\bdm
\frac{3}{4} dT - \frac{1}{2} \sigma^T + \frac{1}{4}\mbox{Scal}^{\nabla} \ = \ \frac{1}{4} dT + \frac{1}{2} \sigma^T + \frac{1}{4}\mbox{Scal}^{\nabla} 
\edm
\bdm
 = \ a \cdot(e_1 \cdot e_2 \cdot e_3 \cdot e_4 + e_1 \cdot e_2 \cdot e_5 \cdot e_6 + e_3 \cdot e_4 \cdot e_5 \cdot e_6 + 3) \
\edm
{\it acts in the spinor modules $\Delta_6^{\pm}$ with eigenvalues} 
$(0, 4 \, a,4 \, a,4\, a)$.\\
 
Consequently, Corollary 3.2 and Theorem 3.3 yield the following results :\\ 

{\bf Theorem 10.8.} {\it Let $(M^6, g, J)$ be a $6$-dimensional nearly K\"ahler spin manifold and let $\nabla$ be the unique linear connection with totally skew-symmetric torsion preserving the nearly K\"ahler structure. Then there exist
 two $\nabla$-parallel spinors. These spinors are sections in the subbundles 
$E^{\pm}$ . If $M^6$ is compact, then every $\nabla$-harmonic spinor is $\nabla$-parallel.}\\

{\bf Remark 10.9.} A $6$-dimensional nearly K\"ahler manifold admits two 
Killing spinors with respect to the Levi-Civita connection (see~\cite{Gru1}) 
and these spinor fields are the $\nabla$-parallel spinors.\\

%\section{References}

\newpage

\small
Thomas Friedrich\\
Humboldt-Universit\"at zu Berlin\\
Institut f\"ur Mathematik\\
Sitz: Rudower Chaussee 25\\
10099 Berlin, Germany\\
{\tt e-mail: friedric@mathematik.hu-berlin.de}\\

Stefan Ivanov\\
University of Sofia ``St. Kl. Ohridski''\\
Faculty of Mathematics and Informatics\\
blvd. James Bourchier 5\\
1164 Sofia, Bulgaria\\
{\tt e-mail: ivanovsp@fmi.uni-sofia.bg}\\

\end{document}